\documentclass[preprint,12pt]{elsarticle}

\usepackage{graphicx}
\usepackage{amssymb}
\usepackage{amsthm}
\usepackage{amsmath}
\usepackage{latexsym}

\newcommand{\abs}[1]{\left|#1\right|}

\journal{}

\begin{document}

\begin{frontmatter}



\title{A variational method for multiphase area-preserving interface motions}


\author[omata,svadlenka2]{Karel Svadlenka\corref{cor1}\fnref{fn}}
\ead{kareru@staff.kanazawa-u.ac.jp}
\author[ginder]{Elliott Ginder}
\ead{eginder@polaris.s.kanazawa-u.ac.jp}
\author[omata]{Seiro Omata}
\ead{omata@se.kanazawa-u.ac.jp}

\cortext[cor1]{Corresponding author}
\fntext[fn]{Tel +81-76-264-5645, Fax +81-76-264-6065}

\address[omata]{Institute of Science and Engineering, Kanazawa University, \\ 
Kakuma-machi, Kanazawa, 920-1192 Japan}
\address[ginder]{Graduate School of Natural Science and Technology, Kanazawa University, \\
Kakuma-machi, Kanazawa, 920-1192 Japan}
\address[svadlenka2]{Institute of Mathematics, Academy of Sciences of the Czech Republic, \\
\v{Z}itn\'{a} 25, 115 67 Praha 1, Czech Republic}

\begin{abstract}
We develop a numerical method for realizing mean curvature motion of interfaces separating multiple phases, whose areas are preserved throughout time.
The foundation of the method is a thresholding algorithm of the Bence-Merriman-Osher type. 
The original algorithm is reformulated in a vector setting, which allows for a natural inclusion of constraints, even in the multiphase case. 
Moreover, a new method for overcoming the inaccuracy of thresholding methods on non-adaptive grids is designed, since this inaccuracy becomes especially prominent in area-preserving motions. 
Formal analysis of the method and numerical tests are presented.
\end{abstract}


\end{frontmatter}


\section{Introduction}

This work develops a method to compute the length-shortening evolution 
of interfaces between an arbitrary number of phases in arbitrary dimension under the constraint that the area of each phase is preserved throughout time.
Such evolutions often appear in situations where interfaces move according to their geometry, while the mass of each phase remains constant  (e.g., grain boundaries in ternary alloys, crystal growth, multiphase flows or formation of soap film bubbles).
In these examples, the motion is driven by the decreasing energy of the internal interfaces, which are out of equilibrium.
This kind of motion also has applications in image processing (denoising, segmentation), in biology (modelling of vesicles and blood cells), in the description of isolated gravitating systems in general relativity \cite{Huisken3}, and other research fields.

Strictly speaking, since area preservation is a global constraint, one cannot consider a constrained curvature flow directly.
Therefore, we have to start from a more basic aspect of the motion, such as the energy.
In particular, we consider the constrained steepest descent of the ``length energy" of each interface, which counts the measure of interfaces weighted by their corresponding interfacial tensions.
The steepest descent of the length energy without any constraint gives the classical mean curvature flow.
On the other hand, in the case of two phases, the area-constrained gradient flow of this energy corresponds to evolution by mean curvature, minus a time-dependent term (equal to the average mean curvature over the interface).
The situation is analogous for more than two phases but the nonlocal term has a complicated form which depends on the configuration of each interface.

The subject is also mathematically interesting, because it is one of the most simple problems with nontrivial limiting behavior. 
It is well known that mean curvature flow shrinks uniformly convex smooth hypersurfaces smoothly to a point in finite time. 
On the other hand, the area-preserving mean curvature flow converges to the solution of the isoperimetric problem, i.e., a sphere \cite{Huisken1, Huisken2, Bellettini, Andrews}.
However, the area-preserving flow may drive general embedded hypersurfaces to self-intersections, as was shown in \cite{Mayer}.
On the other hand, \cite{Escher1} and \cite{Antonopoulou} proved that if the initial surface is sufficiently close to a sphere then it converges to the sphere even if it is not initially convex.
Due to the complexity of the multiphase case, there are only a few results concerning the stability of junctions under area-preserving flow, see \cite{Escher2,Bronsard0} and the references therein.   

Since evolution of surfaces is an intensely studied subject of practical interest, a number of analytical and numerical methods have been developed to treat motion by mean curvature. 
Many of these methods can be applied to the constrained motion addressed here;
let us summarize the known results with emphasis on the multiphase case and volume preservation.

Perhaps the most basic approach is to use the definition of the motion itself. 
That is, in the two-phase case, one computes the evolution of the interface directly from its velocity:
$$ \boldsymbol{v}(x) = ( - \kappa (x) + \kappa_a ) \boldsymbol{n}(x) , \qquad \text{a.e.} \; x \in \partial P(t) , $$
where $P(t)$ denotes the region occupied by one phase, $\kappa$ is the mean curvature and $\kappa_a$ is the average mean curvature over the whole interface.
Algorithms based on this idea are called front tracking methods \cite{Bronsard1,Taylor}. 
They directly approximate the interface based on the Huygens' principle and are effective for computing the evolution of smooth surfaces without topological changes. 
Although this method is simple in principle, if interaction of different parts of the interface occurs, a complicated decision algorithm is necessary to proceed with the computation, and this becomes increasingly more involved in higher dimensions. Higher dimensions and a larger number of phases also complicate the calculation of curvatures and their averages over the interfaces.

A more general framework is provided by the level set approach which, thanks to its implicit representation of the interface, is able to deal with topological singularities and nonsmooth data.
In this approach, the initial interface $\partial P(0)$ is expressed as the $0$-level set of a function $\phi(x,0)$, and the mean curvature flow is achieved as 
$$ \partial P(t) = \{ x; \phi(x,t)=0  \} , $$
where $\phi$ is the viscosity solution of the Hamilton-Jacobi equation
$$ \phi_t = \text{div} \; \Big( \frac{\nabla \phi}{| \nabla \phi |} \Big) | \nabla \phi | . $$
The constrained flow can be realized in this setting by considering the area-constrained gradient flow of the length energy functional written in terms of the level set function
$$ L(\partial P(t)) = \int \delta(\phi(x,t)) | \nabla \phi (x,t)| \, dx, $$
where $\delta$ denotes the Dirac delta function. 
It is still necessary to calculate the curvature values, but this method can be extended to the multiphase setting by introducing as many level set functions as there are phases and imposing an additional constraint so that the level sets do not overlap or create vacuums (see \cite{Zhao} or \cite{Zhao2}).
However, such a constraint has an unwanted impact on the flow \cite{Esedoglu} and the phase areas are not adequately preserved during the computations. 
The first problem was solved in \cite{Esedoglu} and \cite{Kublik} by employing signed distance functions. 
Area-preserving motions are also addressed in these works but are limited to two-phase flows. 
Another multiphase modification of the level-set method, which has some similarities to our own approach, was developed in \cite{Smith}, where the constraint of \cite{Zhao} was replaced by a projection step. 
However, the impact of the projection on the dynamics of the interface was not analyzed.

The area-preserving mean curvature flow arises as a limit of the following nonlocal mass-preserving diffusion equation \cite{Golovaty, Bronsard2, Brassel}
$$ u_t = \Delta u - \frac{1}{\varepsilon^2} W'(u) + \frac{1}{\varepsilon^2 | \Omega| } \int_{\Omega} W'(u) \, dx , $$
where $W$ is a double-well potential and $\varepsilon$ is a small parameter related to the width of the diffuse interface.
It has been shown that, under suitable conditions, the set 
$$ P_{\varepsilon}(t) = \{ x; \; u_{\varepsilon}(x,t) \geq \tfrac{1}{2} \} $$
approximates $P(t)$ with error $O(\varepsilon^2 |\log \varepsilon|^2)$.

Based on this fact, the so-called phase field methods  represent interfaces by thin layers in the solution and thus the resolution of this internal layer requires a very fine mesh. 
On the other hand, this approach handles topological changes without trouble 
and does not require explicit computation of curvatures. 
An interesting computational analysis for the multiphase case is given in \cite{Garcke}.


The basic idea of this paper is to use another approach, often referred to as a thresholding method. We adapt the so-called BMO algorithm (sometimes also called the MBO algorithm), which was discovered by Merriman, Bence and Osher in \cite{MBO}, to generate multiphase area-preserving motion.
As far as is known to the authors, the existing works (except \cite{Zhao2}) do not treat the approximation of multiphase length-shortening flow under area constraints.

The BMO algorithm exploits the fact that short-time diffusion of the characteristic function of a region enclosed by an interface (i.e., its convolution with the Gaussian kernel), evolves the interface according to its mean curvature.
More precisely, the characteristic function of a region is evolved for a short time by the heat equation and then a thresholding step is carried out to obtain the new interface (given by the $1/2$-level set of the diffused function).
The main advantage of this approach is that it naturally treats topological changes, produces no intercalary regions and does not require explicit computation of curvatures.
Moreover, it is numerically attractive because of its stability and low computational complexity.

This thresholding method was applied to multiphase flow in \cite{MBO}, while \cite{Ruuth} uses the BMO algorithm for constructing two-phase area-preserving curvature flow. 
However, the latter method lacks theoretical support, and it is not clear how to extend it to more than two phases and to more general motions. 
We therefore introduce a different approach, which is also based on the BMO method. Our method can treat any number of phases in any dimension and can be extended to more general motions, such as mean curvature motion with transport.

The difficulty of the multiphase constrained flow is that the phases influence each other not only locally via the shape of their interface, but also globally via their areas.
The idea used to overcome this complication is to formulate the original multiphase BMO method in a vector-valued fashion and to realize the area constraint by considering a constrained gradient flow.
This constrained flow presents a computational difficulty, due to the fact that the interfacial velocities are slower when compared to the flow without area preservation.
That is, since the interface must move at least the distance of the grid size at each time step, this places unreasonable restrictions on the grid resolution used in the numerical implementation.
We are able to overcome these restrictions by introducing a technique of temporary and localized refinement.

The paper is organized in the following way. 
In section 2, we introduce the area-preserving mean curvature flow, as well as the BMO algorithm for its approximation. 
We discuss the numerical algorithm in section 3, and section 4 concerns its implementation. 
Section 5 presents a number of numerical examples and analyses of errors and model parameters.
The appendix includes a number of theoretical results requiring technical computations.

\section{BMO algorithm for area-preserving mean curvature flow}

\subsection{Area-preserving mean curvature flow}

\subsubsection{Two-phase case}
\label{mcf_two}

We first fix the meaning of the terms ``area" and ``length" used in the sequel.
Working in $m$-dimensional space, the word ``area" shall mean the $m$-dimensional Lebesgue measure of the region corresponding to a phase, i.e., in 3 dimensions it regards the volume. 
The word ``length", on the other hand, shall always refer to the $(m-1)$-dimensional Lebesgue measure of the boundary of a phase region, thus in 3 dimensions corresponding to the surface area of the interface.

\begin{figure}[!ht]
\begin{center}
\label{fig_2phases}
\includegraphics[scale=0.7]{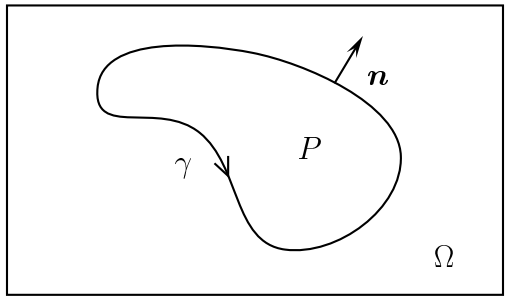}
\caption{Two phases divided by an interface.}
\end{center}
\end{figure}

Mean curvature flow is related to systems whose energy depends on the length of their surface, e.g., soap films. 
We will explain the basic equations for curves in the two-dimensional plane, since the derivations in higher dimensions are similar but lack transparency. 
Accordingly, let us consider a smooth Jordan curve $\gamma$ contained in a subregion $\Omega$ as in figure 1.
Let the curve be parametrized:
$$ \gamma (s) = (\gamma_1(s), \gamma_2(s)), \quad s \in [a,b], \qquad \gamma(a)=\gamma(b) , $$
and in such a way that the enclosed region $P$ is on the left side of the curve.
Then the length of the curve is given by
$$ L(\gamma) = \int_a^b | \gamma'(s) | \, ds . $$
The gradient flow of the above energy can be found from its first variation. 
That is, for any smooth closed curve $\varphi(s), s \in [a,b]$ we compute
\begin{equation}
\label{varlen} 
\frac{d}{d \varepsilon} L( \gamma + \varepsilon \varphi)|_{\varepsilon = 0} = \int_{\gamma} (\kappa \boldsymbol{n}) \cdot \varphi \, dl , 
\end{equation}
where $\kappa$ is the curvature, $\boldsymbol{n}$ is the unit outer normal to the curve $\gamma$ at a given point and where we integrate with respect to arc length $l$:
$$ \kappa (s) = \frac{\gamma_1'\gamma_2''-\gamma_2'\gamma_1''}{|\gamma'|^3}(s), \qquad \boldsymbol{n}(s) = \frac{1}{|\gamma'|} (\gamma_2', - \gamma_1')(s), \qquad dl= |\gamma'(s)| ds. $$
In the general $m$-dimensional case one obtains the same formula as (\ref{varlen}), where $\kappa$ means the trace of second fundamental form divided by $m$.
Hence, the fastest shortening of the curve occurs for the flow with normal velocity equal to minus (a multiple of) the mean curvature.

Next we consider the curve-shortening flow under the constraint of area preservation.
The area functional for region $P$ reads
$$ A(\gamma) = \frac{1}{2} \int_{\gamma} (x,y) \cdot \boldsymbol{n} \, dl = \frac{1}{2} \int_a^b (\gamma_1 \gamma_2'- \gamma_2 \gamma_1') \, ds $$
and its first variation is
$$ \frac{d}{d \varepsilon} A(\gamma + \varepsilon \varphi) |_{\varepsilon = 0} = \int_{\gamma} \boldsymbol{n} \cdot \varphi \, dl .$$
Analogously, in any dimension the variation turns out to be the normal vector at each point of the hypersurface.

Following the construction in \cite{Svadlenka}, we introduce a time-dependent Lagrange multiplier $\lambda(t)$ for the area constraint and express the velocity of the area-constrained mean curvature flow by
$$ \boldsymbol{v} =  (-\kappa + \lambda) \boldsymbol{n} . $$
A precise expression for the the multiplier $\lambda$ can be obtained in the following standard way.
From the fact that the area is preserved one has
$$ \frac{d}{dt} A(\gamma(t)) = \int_{\gamma(t)} \boldsymbol{v}(t) \cdot \boldsymbol{n}(t) \, dl = 0 .$$
Hence it follows that the integral of $-\kappa + \lambda$ over the curve $\gamma$ vanishes at each time. This yields
$$ \lambda (t) = \frac{1}{L(\gamma(t))} \int_{\gamma(t)} \kappa (t) \, dl. $$
That is, the Lagrange multiplier expresses the average mean curvature along the interface.

\subsubsection{Multiphase case}

We briefly derive the velocity of interfaces moving by area-preserving mean curvature flow in the multiphase two-dimensional setting. 
We assume that the number of phases $k$ is finite and that the interfaces between different pairs of phases form a finite collection of arcs.
The boundary $\gamma_i$ of phase region $P_i$ can then be written as the union of the interfaces from all other phase regions:
$$ \gamma_i = \bigcup_{j \neq i} \gamma_{ij} = \bigcup_{j \neq i} \bigcup_l \gamma_{ij}^l, \qquad i=1, \dots, k .$$
Here $\gamma_{ij}$ denotes the interface between phases $P_i$ and $P_j$, and the index $l$ expresses the fact that each interface may have several disjoint parts.
In the following, however, we omit the decomposition using index $l$ since it has no influence on the computations.

Each interface $\gamma_{ij}^l$ is considered as an oriented curve, which has the region $P_i$ on its left side.
For any point interior to the interface $\gamma_{ij}$ we define the normal $\boldsymbol{n}$ as the unit vector pointing into the phase with larger index, i.e., $\boldsymbol{n}$ is the outer normal to $P_i$ if $i<j$.
We remark that some of the curves $\gamma_{ij}$ may be empty.

The energy of the system, considering arbitrary surface tension $\tau_{ij}$ for $\gamma_{ij}$, is
\begin{equation}
\label{lemul}
\sum_{i=1}^k \sum_{j>i} \tau_{ij} \int_{\gamma_{ij}} \, dl . 
\end{equation}
This value is to be minimized under the condition of constant areas, that is,
$$ \frac{1}{2} \Big( \sum_{j>i} \int_{\gamma_{ij}} \boldsymbol{x} \cdot \boldsymbol{n} \, dl - \sum_{j<i} \int_{\gamma_{ji}} \boldsymbol{x} \cdot \boldsymbol{n} \, dl \Big) = \text{const} , \qquad i = 1, \dots , k-1 . $$
Introducing Lagrange multipliers $\lambda_i$, $i=1, \dots ,k-1$, for each of the constraints, the constrained variation reads
$$ \sum_{i=1}^k \sum_{j>i} \tau_{ij} \int_{\gamma_{ij}} \kappa \boldsymbol{n} \cdot \varphi^{ij} - \sum_{i=1}^{k-1} \lambda_i \left( \sum_{j>i} \int_{\gamma_{ij}} \boldsymbol{n} \cdot \varphi^{ij} \, dl  - \sum_{j<i} \int_{\gamma_{ji}} \boldsymbol{n} \cdot \varphi^{ij} \, dl \right), $$
where $\varphi^{ij}, \, i,j=1, \dots ,k$, denotes a smooth perturbation within $\gamma_{ij}$.

From this it follows that the magnitude of normal velocity $v^{ij}$ for interface $\gamma_{ij}$ ($i<j$) in the direction of $\boldsymbol{n}$ will be
$$ v^{ij}=-\tau_{ij} \kappa + \lambda_i - (1-\delta_{jk}) \lambda_j . $$
Here $\delta_{jk}$ denotes the Kronecker delta, which arises from the redundance of the area constraint for phase $P_k$ (it follows from the constraints on the other phases and the fixation of the domain).

Explicit representations of the Lagrange multipliers can be obtained as follows.
Let us denote by $L_{ij}$ the length of interface $\gamma_{ij}$, by $L_{i}$ the length of the boundary of phase region $P_{i}$, and by $\kappa_{ij}$ the average mean curvature through interface $\gamma_{ij}$ times the tension $\tau_{ij}$:
$$ \kappa_{ij} = \frac{\tau_{ij}}{L_{ij}} \int_{\gamma_{ij}} \kappa \, dl . $$ 
Then the preservation of phase area $P_i$ gives the condition
\begin{eqnarray*}
0 &=& \sum_{j>i} \int_{\gamma_{ij}} v^{ij} \, dl - \sum_{j<i} \int_{\gamma_{ji}} v^{ji} \, dl \\
&=& \sum_{j>i} \int_{\gamma_{ij}} (-\tau_{ij} \kappa ) \, dl - \sum_{j<i} \int_{\gamma_{ji}} (-\tau_{ij} \kappa ) \, dl \\
&& + \sum_{j>i} \lambda_i \int_{\gamma_{ij}} \, dl + \sum_{j<i} (1-\delta_{ik}) \lambda_i \int_{\gamma_{ji}} \, dl \\
&& + \sum_{j>i}  (1-\delta_{jk}) (-\lambda_j) \int_{\gamma_{ij}} \, dl - \sum_{j<i} \lambda_j \int_{\gamma_{ji}} \, dl \\
&=& - \sum_{j>i} L_{ij} \kappa_{ij} + \sum_{j<i} L_{ij} \kappa_{ij} + \lambda_i \sum_{j \neq i} L_{ij} - \sum_{j \neq i} (1-\delta_{jk}) \lambda_j L_{ij} . 
\end{eqnarray*}
Since the above holds for $i=1, \dots, k-1$, we obtain a system of linear equations for $\lambda_1, \dots , \lambda_{k-1}$:
$$ L_{i} \lambda_i - \sum_{\substack{j=1\\j \neq i} 
}^{k-1} L_{ij} \lambda_j  =  \sum_{j>i} L_{ij} \kappa_{ij} - \sum_{j<i} L_{ij} \kappa_{ij}  , \quad i= 1, \dots , k-1. $$ 

For a given configuration, the solution to this system gives the Lagrange multipliers. For example, in the case of  three phases, as in figure 2, one has the velocities:
\begin{eqnarray*}
v^{13} &=& -\tau_{13} \kappa + \big( 1- \tfrac{L_{13} L_{12}}{\alpha} \big) \kappa_a^{P_1} - \big( 1- \tfrac{L_{13}L_{23}}{\alpha} \big) \kappa_a^{P_3}  , \\
v^{23} &=& -\tau_{23} \kappa + \big( 1- \tfrac{L_{23}L_{12}}{\alpha} \big) \kappa_a^{P_2} - \big( 1- \tfrac{L_{13}L_{23}}{\alpha} \big) \kappa_a^{P_3}  , \\
v^{12} &=& - \tau_{12} \kappa + \big( 1- \tfrac{L_{13}L_{12}}{\alpha} \big) \kappa_a^{P_1} - \big( 1- \tfrac{L_{23}L_{12}}{\alpha} \big) \kappa_a^{P_2} . 
\end{eqnarray*}
Here $\alpha = L_{13}L_{23} + L_{23}L_{12} + L_{13}L_{12}$ and $\kappa_a^{P_i}$, $i=1,2,3$, represent the average mean curvatures along the whole boundary of each phase, weighted by the surface tension.

\begin{figure}[!ht]
\begin{center}
\label{fig_3phases}
\includegraphics[scale=0.7]{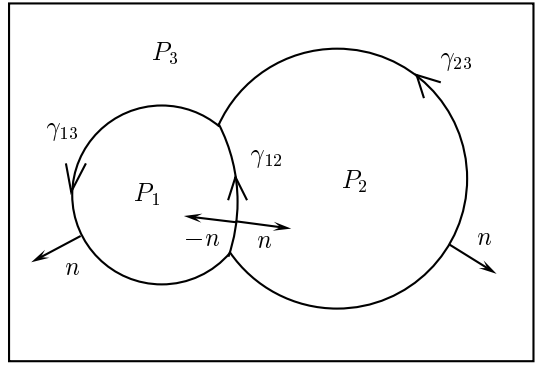}
\caption{An example of a three-phase configuration.}
\end{center}
\end{figure}

When phases $P_1$ and $P_2$ are separated, i.e., when $L_{12}=0$, the above formulae reduce to the form of the two-phase flow from section \ref{mcf_two}.
The extension of the above calculations to higher space dimensions is natural, as we use only the notions of interfaces and their oriented normals.

\subsection{BMO algorithm}

For the sake of clarity, we explain our method in three successive steps.
First we summarize the original idea of the BMO algorithm for two phases when no area constraint is present.
Then we describe the existing algorithm for an arbitrary number of phases, again without area constraint, and reformulate the algorithm in a vector-type setting.
In the third step we finally design the method for multiphase area-preserving motion.

\subsubsection{Two-phase motion without area constraint}

We describe the BMO algorithm for the case when only two phases are present. 
This algorithm works in any space dimension. Given an initial interface $\gamma$, we take $P$ to be the region enclosed by this interface (and possibly by the boundary of the region $\Omega$ where the motion is considered) and define its characteristic function $\chi$ as 
$$ \chi(x) = \left\{ \begin{array}{l} 1 \qquad x \in P, \\ 0 \qquad x \not\in P. \end{array} \right. $$
The new interface after a time $\Delta t$ will be the boundary of the $\tfrac{1}{2}$-level set of the solution to the heat equation at time $\Delta t$ with initial datum $\chi$. 
The two-phase algorithm thus reads as follows:
\begin{enumerate}
  \item Given a region $P$, set $\chi$ to be its characteristic function.
  \item \label{2step2} Solve the heat equation with initial condition $\chi$:
     \begin{eqnarray*}
     u_t(t,x) &=& \Delta u(t,x) \qquad \;\, \text{for} \; (t,x) \in (0, \Delta t] \times \Omega , \\
     \frac{\partial u}{\partial n} (t,x) &=& 0 \qquad \qquad \quad \;\, \text{on} \; (0, \Delta t] \times \partial \Omega , \\ 
     u(0,x) &=& \chi(x) \qquad \qquad \text{in} \; \Omega .
     \end{eqnarray*}
  \item Update $\chi$ as the $\tfrac{1}{2}$-level set of $u(\Delta t, \cdot)$:
     $$ \chi (x) = \left\{ \begin{array}{l} 1 \qquad \text{if} \;\; u(\Delta t, x) > \tfrac{1}{2}, \\ 0 \qquad \text{if} \;\; u(\Delta t,x) \leq \tfrac{1}{2} . \end{array} \right.  $$
The evolved interface is now the boundary of the set $\{x \in \Omega ; \; \chi (x) =1 \}$.
  \item Go back to step \ref{2step2} to proceed with the computation for the next time step.
\end{enumerate}

It has been rigorously shown, in a general framework including topological changes, that this algorithm converges to motion by mean curvature as $\Delta t \to 0$ \cite{Barles, Evans, Goto}.

Here we remark that the Neumann boundary condition in the diffusion step guarantees that the interface will touch the boundary of $\Omega$ with right angle. Other boundary conditions, such as Dirichlet conditions pinning the interface at the boundary, may be used according to necessity. \\

\subsubsection{Multiphase motion without area constraint}
\label{multinoBMO}

We next address the case of multiple phases.
The idea of sharpening separately diffused characteristic functions (one for each phase) was introduced in \cite{MBO}. 
The algorithm is as follows.

\begin{enumerate}
  \item Given regions $P_i, i=1, \dots ,k$, set $\chi_i$ to be the characteristic function of $P_i$.
  \item \label{kstep2} For $i=1, \dots ,k$, obtain $u_i(\Delta t,x)$ by solving the heat equation with initial condition $\chi_i$ up to time $\Delta t$. 
  \item Update $\chi_j$ as the characteristic function of the set where $u_j$ has the largest value amongst all solutions  $u_i$:
     $$ \chi _j(x) = \left\{ \begin{array}{l} 1 \qquad \text{if} \;\; u_j(\Delta t, x) \geq u_i(\Delta t,x) \;\; \forall i \ne j, \\ 0 \qquad \text{otherwise} . \end{array} \right.  $$
The new interfaces are the boundaries of the sets $\{x \in \Omega ; \; \chi_i (x) =1 \}$.
  \item Go back to step \ref{kstep2} to proceed with the computation for the next time step.
\end{enumerate}

The above algorithm can be reformulated to obtain an equivalent algorithm using a single vector-valued heat equation. 
This is essential for implementing the area constraint and for dealing with more general motions. 

We prepare $k$ reference unit vectors $\boldsymbol{p}_i, i= 1, \dots, k$, of dimension $k-1$, each corresponding to a phase $P_i$. 
They are defined as the vectors pointing from the centroid of a standard $k$-simplex to its vertices (cf. figure A.12 in the appendix). 
Hence, there are infinitely many possible $k$-tuples but the relative distributions of the vectors are identical. 
See \ref{app_refvec} for a simple way to construct these vectors.

Using the vectors $\boldsymbol{p}_i$, the multiphase algorithm can be written as follows:
\begin{enumerate}
  \item Given regions $P_i, i = 1, \dots ,k$, set
  			$$ 	\boldsymbol{u}_0(x) = \boldsymbol{p}_i \qquad \text{for} \;\; x \in P_i . $$
  \item \label{knstep2} Solve the vector-valued heat equation with initial condition $\boldsymbol{u}_0$:
     \begin{eqnarray}
		 \label{vech}
     \boldsymbol{u}_t(t,x) &=& \Delta \boldsymbol{u}(t,x) \qquad \;\; \text{for} \; (t,x) \in (0, \Delta t] \times \Omega , \\
     \frac{\partial \boldsymbol{u}}{\partial n} (t,x) &=& 0 \qquad \qquad \quad \;\;\; \text{on} \; (0, \Delta t] \times \partial \Omega , \nonumber \\ 
     \boldsymbol{u}(0,x) &=& \boldsymbol{u}_0(x) \qquad \qquad \text{in} \; \Omega . \nonumber
     \end{eqnarray}
  \item Update $\boldsymbol{u}_0$ by identifying the reference vector which is closest to the solution $\boldsymbol{u}(\Delta t, x)$:
		  	\begin{equation}\label{truncstep}
     \boldsymbol{u}_0 (x) = \boldsymbol{p}_j, \qquad \text{where} \quad \boldsymbol{p}_j \cdot \boldsymbol{u}(\Delta t, x) = \max_{i=1, \dots , k} \boldsymbol{p}_i \cdot \boldsymbol{u}(\Delta t, x) .
     			\end{equation}
This redistribution of reference vectors determines the approximate new phase regions after time $\Delta t$.
  \item Go back to step \ref{knstep2} to proceed with the computation for the next time step.
\end{enumerate}

The equivalence of our algorithm with the original BMO can be shown by considering the functions 
$$ w_i (t,x) = \frac{k-1}{k} \left( \boldsymbol{u}(t,x) \cdot \boldsymbol{p}_i + \frac{1}{k-1} \right) , \quad i=1, \dots ,k. $$
Since $\boldsymbol{u}$ is the solution to the linear vector-valued heat equation, we have
\begin{eqnarray*}
(w_i)_t &=& \Delta w_i \qquad \text{in} \; (0,T] \times \Omega , \\
\frac{\partial w_i}{\partial n} &=& 0 \qquad \quad \; \text{on} \; (0,T] \times \partial \Omega .
\end{eqnarray*}
Moreover, since from the construction of $\boldsymbol{p}_i$ it holds that $\boldsymbol{p}_i \cdot \boldsymbol{p}_j = 1/(1-k)$ for $i \neq j$, we can readily check the identity
$$ w_i(0,x) = \frac{k-1}{k} \left( \boldsymbol{u}_0(x) \cdot \boldsymbol{p}_i + \frac{1}{k-1} \right) = \chi_i(x), $$
where $\chi_i$ is the characteristic function of $i$-th phase region.
It follows that $w_i$ is identical to the solution $u_i$ of the scalar heat equation from the original algorithm for each $i=1, \dots ,k$.
From the definition of $w_i$ it is immediate that
$$ \boldsymbol{u}(x,t) \cdot \boldsymbol{p}_i \geq \boldsymbol{u} (x,t) \cdot \boldsymbol{p}_j \quad \Leftrightarrow \quad w_i (x,t) \geq w_j (x,t) \quad \Leftrightarrow \quad u_i(x,t) \geq u_j(x,t), $$
which proves the equivalence of both algorithms.  \\

\noindent
{\bf Remark.} The reference vectors are related to the position of wells in the phase field approach.
Indeed, the idea of the BMO algorithm for the case of two phases originated in a simple splitting scheme for the singularly perturbed reaction-diffusion equation
$$ u_t = \varepsilon \Delta u - \tfrac{1}{\varepsilon} W'(u), $$
where $W$ is a double-well potential.
Here, the splitting scheme consists of two steps. 
The first step solves the heat equation $u_t= \varepsilon \Delta u$, which corresponds to the diffusion step of the BMO algorithm, and the second solves $u_t = - \frac{1}{\varepsilon} W'(u)$. 
This corresponds to the thresholding step if the equation is solved for sufficiently long time.
Here we can see that the thresholding values $0$ and $1$ in the BMO algorithm correspond to the positions of the two wells of $W$.
Accordingly, if we want to calculate three-phase motion, we can look at the problem from the viewpoint of constructing a suitable well-type potential. 
A potential with three wells at different positions along a scalar axis would obviously yield incorrect results, because the strength of the wells would not be equivalent.
Therefore, we have to increase the number of variables for the potential and construct the wells in a symmetric way. 
The reference vectors introduced above then give the coordinates of the positions of the wells.

\subsubsection{Multiphase motion with area constraint}
\label{mmac}

The paper \cite{Ruuth} presents an area-preserving BMO method for 2 phases.
The authors adjust the height at which the thresholding occurs in such a way that the resulting area of the level set is preserved. 
Such a level set is guaranteed to exist by the maximum principle and one can compute that the thresholding height has to be changed from $\tfrac{1}{2}$ to
$$ \frac{1}{2} - \frac{1}{2} \kappa_a \sqrt{\frac{\Delta t}{\pi}} , $$
where $\kappa_a$ is the average mean curvature along the interface.

However, when three or more phases are present, each interface has a different average mean curvature and the thresholding heights become different. 
One could try to use a different thresholding height for each phase but then the global interaction would be ignored and the phases may overlap or create vacuums, especially when they are initially touching.
Therefore, we suggest a different approach incorporating the area constraint into the diffusion process.
In this method, a heat equation with a nonlinear source term expressing the area preservation of level sets is solved. Subsequently, the solution is sharpened at a fixed height.  

The mentioned nonlinear heat equation corresponds to the gradient flow of the functional 
$$ J(\boldsymbol{u}) = \int_{\Omega}  | \nabla \boldsymbol{u} |^2 \, dx $$
in the constrained set
$$ {\mathcal K} = \Big\{ \boldsymbol{u} \in H^1(\Omega; {\bf R}^{k-1});  \\
\int_{\Omega} \chi_{ \{ \boldsymbol{u}(x) \cdot \boldsymbol{p}_i \geq \boldsymbol{u}(x) \cdot \boldsymbol{p}_j \; \forall j \}} \, dx = A_i \;\; \text{for} \; i=1, \dots , k-1 \Big\} , $$
where $A_i$ denotes the given area of phase $P_i$.
For simplicity, we consider only the case when one type of interface divides one type of phase into multiple regions, like in a soap froth. 
When multiple phases and interfaces with nonuniform properties are present, a different form of the energy $J$ has to be adopted \cite{Malek}.

Let us consider the two-dimensional two-phase case as in figure 1 to understand the meaning of this gradient flow. 
Our argument is formal and we simplify the exposition by discretizing time and adopting the numerical approach from the next section, which is based on the method of Rothe \cite{Rothe}. 

Let the curve $\gamma$ denote the $\tfrac{1}{2}$-level set of a scalar function $u$ and call the region enclosed by this curve $P$. 
We assume that an initial function $u_*$ is given and we search for the minimizer of the integral
\begin{equation}
\label{mincon}
J(u) = \int_{\Omega} \left( \frac{|u-u_*|^2}{h} + | \nabla u|^2 \right) \, dx 
\end{equation}
under the constraints
\begin{eqnarray}
\label{c1} u(\gamma(s)) &=& \frac{1}{2}, \qquad s \in [a,b], \\
\label{c2} \text{meas}(P) &=& A .
\end{eqnarray}
Here $h$ denotes the length of the discrete time step and $A$ is the required area.

Note that by perturbing function $u$ in regions away from the $\frac{1}{2}$-level set, one can readily deduce that the minimizer $u$ satisfies, in a weak sense , the following
\begin{eqnarray}
\label{1varJ}
\frac{u-u_*}{h} - \Delta u &=& 0 \qquad \text{in} \; P \cup (\Omega \setminus \bar{P}), \\
\frac{\partial u}{\partial n} &=& 0 \qquad \text{on} \; \partial \Omega . \nonumber
\end{eqnarray}

Now consider a perturbation of $u$ of the form $u+ \delta u$, which is allowed to affect the $\frac{1}{2}$-level set.
The corresponding change in the level curve can be written in the form $\alpha (\gamma + \delta \gamma)$, $\delta \gamma = (\delta \gamma_1, \delta \gamma_2)$, where $\alpha$ is a constant depending on $\delta \gamma$  and whose role is to adjust the area to the correct value. 
Because of constraints (\ref{c1}) and (\ref{c2}), we cannot choose $\delta u$ arbitrarily.
However, we can select an arbitrary $\delta \gamma$, find an appropriate constant $\alpha$, and use the corresponding $\delta u$.

We compute the variation of functional $J$ remembering that $u$ may not be smooth across the interface $\gamma$:
\begin{eqnarray*}
&& J(u+ \delta u) - J(u) \\
&& \qquad \simeq \int_{\Omega} \left( \frac{u-u_*}{h} \delta u + \nabla u \nabla (\delta u) \right) \, dx \\
&& \qquad \simeq \int_{P} \left( \frac{u-u_*}{h} \delta u + \nabla u \nabla \delta u \right) \, dx + \int_{\Omega \setminus \bar{P}} \left( \frac{u-u_*}{h} \delta u + \nabla u \nabla \delta u \right) \, dx \\
&& \qquad \simeq \int_{P} \left( \frac{u-u_*}{h} - \Delta u \right) \delta u  \, dx + \int_{\gamma} \frac{\partial u_P}{\partial n} \delta u_P \, dl \\
&& \qquad \; + \int_{\Omega \setminus \bar{P}} \left( \frac{u-u_*}{h} - \Delta u \right) \delta u  \, dx - \int_{\gamma} \frac{\partial u_{\Omega \setminus \bar{P}}}{\partial n} \delta u_{\Omega \setminus \bar{P}} \, dl + \int_{\partial \Omega} \frac{\partial u}{\partial n} \delta u \, dl . 
\end{eqnarray*}
Here $w_P$ denotes the value of a function $w$ taken as a limit from inside of $P$. 
The symbol $w_{\Omega \setminus \bar{P}}$ has the analogous meaning.

Due to (\ref{1varJ}), on the interface we have
$$ \int_{\gamma} \left( \frac{\partial u_P}{\partial n} \delta u_P - \frac{\partial u_{\Omega \setminus \bar{P}}}{\partial n} \delta u_{\Omega \setminus \bar{P}} \right) \, dl = 0 \qquad \text{for all admissible} \; \delta u . $$
We next to rewrite this identity in terms of the arbitrary perturbation $\delta \gamma$.
To this end, we use the fact that the phase area is preserved to obtain
$$ \frac{1}{2} \alpha^2 \int_a^b \left[ (\gamma_1+ \delta \gamma_1)(\gamma_2'+\delta \gamma_2') - (\gamma_2+\delta \gamma_2)(\gamma_1'+ \delta \gamma_1') \right] \, ds = A, $$ 
which yields
$$ \int_{\gamma} \delta \gamma \cdot \boldsymbol{n} \, dl \simeq \frac{1- \alpha^2}{\alpha^2} A . $$
Here, the symbol $\simeq$ means that the equation holds up to second order in terms of $\delta \gamma$ or $\delta u$.
From this we can compute the value of $\alpha -1$ which will be needed later:
\begin{equation}
\label{am1}
\alpha -1 \simeq - \frac{1}{2A} \int_{\gamma} \delta \gamma \cdot \boldsymbol{n} \, dl . 
\end{equation}

Condition (\ref{c1}) means
$$ (u+ \delta u)(\alpha (\gamma + \delta \gamma)) - u(\gamma) = 0 .$$
Using a Taylor expansion and (\ref{am1}) we obtain that, on $\gamma$, one has
\begin{eqnarray}
0 & \simeq & \nabla u_P \cdot \big( (\alpha -1) \gamma + \delta \gamma \big) + \delta u_P \nonumber \\
\label{du1} & \simeq & \nabla u_P \cdot \Big( \delta \gamma - \frac{\gamma}{2A} \int_{\gamma} \delta \gamma \cdot \boldsymbol{n} \, dl \Big) + \delta u_P \\
& \simeq & \nabla u_{\Omega \setminus \bar{P}} \cdot \Big( \delta \gamma - \frac{\gamma}{2A} \int_{\gamma} \delta \gamma \cdot \boldsymbol{n} \, dl \Big) + \delta u_{\Omega \setminus \bar{P}}. \nonumber
\end{eqnarray}
Expressing $\delta u$ from (\ref{du1}) as
\begin{eqnarray*}
\delta u_P &=& - \frac{\partial u_{P}}{\partial n} \boldsymbol{n} \cdot \left( \delta \gamma - \frac{\gamma}{2A} \int_{\gamma} \delta \gamma \cdot \boldsymbol{n} \, dl \right) , \\
\delta u_{\Omega \setminus \bar{P}} &=& - \frac{\partial u_{\Omega \setminus \bar{P}}}{\partial n} \boldsymbol{n} \cdot \left( \delta \gamma - \frac{\gamma}{2A} \int_{\gamma} \delta \gamma \cdot \boldsymbol{n} \, dl \right), 
\end{eqnarray*}
we obtain
$$ \int_{\gamma} \left[ - \big( \frac{\partial u_{P}}{\partial n} \big)^2 + \big( \frac{\partial u_{\Omega \setminus \bar{P}}}{\partial n} \big)^2 \right] \boldsymbol{n} \cdot \left( \delta \gamma - \frac{\gamma}{2A} \int_{\gamma} \delta \gamma \cdot \boldsymbol{n} \, dl \right) \, dl = 0
\qquad \forall \delta \gamma . $$
Noting that 
$$ \int_{\gamma} \boldsymbol{n} \cdot \left( \delta \gamma - \frac{\gamma}{2A} \int_{\gamma} \delta \gamma \cdot \boldsymbol{n} \, dl \right) \, dl = 0 ,$$
we arrive at the interface condition
$$ \left( \frac{\partial u_{P}}{\partial n} \right)^2 -  \left( \frac{\partial u_{\Omega \setminus \bar{P}}}{\partial n} \right)^2  = \lambda = \text{const} \qquad \text{on} \; \gamma . $$
Similar calculation can be carried out in the multiphase setting, see \ref{app_multi} for more details.

In view of the derived condition on $\gamma$ and the results regarding two-phase free boundary problems in \cite{Alt}, we can reformulate the variational problem in the following way:
Find a minimizer $u_{\lambda} \in H^1(\Omega)$ of the functional 
\begin{equation}
\label{minlam}
J^{\lambda}(u) = \int_{\Omega} \left( \frac{|u-u_*|^2}{h} + | \nabla u |^2 + \lambda \chi_{u> \frac{1}{2}} \right) \, dx . 
\end{equation}
According to the results in \cite{Alt}, we can expect that this minimization problem is in a sense equivalent to looking for a weak solution of the problem
\begin{eqnarray}
\label{nonlinheat}
u_t - \Delta u &=& \mu \qquad \text{in} \; (0,T] \times \Omega , \\
\frac{\partial u}{\partial n} &=& 0 \qquad \text{on} \; \partial \Omega , \nonumber
\end{eqnarray}
where $\mu$ is a Radon measure given by
$$ \mu (t,x) = \lambda (t) {\mathcal H}^1 \lfloor_{\partial P(t)}, \qquad P(t) = \{x; \; u(t,x) > \tfrac{1}{2} \}, $$
for a suitable space-independent function $\lambda$.
Here the symbol ${\mathcal H}^1$ means the one-dimensional Hausdorff measure.  
This type of problem is known as the two-phase parabolic free boundary problem \cite{Caffarelli1, Caffarelli2, Weiss}. 

Having deduced the partial differential equation corresponding to the constrained gradient flow, let us return to the considerations concerning the BMO method.
Formal calculation shows that the BMO method using the parabolic equation (\ref{nonlinheat}) instead of the heat equation evolves interfaces with normal velocity equal to minus mean curvature, plus a space-independent term:
$$ v = -\kappa + 2 \lambda (0) + O(t), \qquad t \to 0+.$$ 
Thus, if $\lambda$ is chosen appropriately, the area-preserving mean curvature flow is realized.
Since the derivation of the above relation is rather technical, we present it in \ref{appA}.

Minimization of (\ref{minlam}) is much easier to implement than the constrained minimization (\ref{mincon}) or the partial differential equation (\ref{nonlinheat}).
In each time step of the BMO algorithm, if for some $\lambda$ we obtain that the area of $\{x; \; u_{\lambda}(x) > \tfrac{1}{2} \}$ is equal to the given value $A$, then $u_{\lambda}$ becomes a solution to our original problem.
However, it is not clear how to calculate such a $\lambda$ or even if such a function exists.
For the two-phase case we see that $\lambda$ should be half of the average mean curvature, but determining the Lagrange multiplier for the multiphase case is complicated.
Moreover, one of the advantages of the BMO approach is that it does not require the computation of curvature values, so we want to avoid the direct calculation of $\lambda$.

Therefore, based on the results in \cite{Aguilera} we approach (\ref{mincon}) by using the method of penalization and consider minimization of the functional
$$ J_{\varepsilon} (u) = \int_{\Omega} \Big( \frac{|u-u_*|^2}{h} + | \nabla u |^2 \Big) \, dx + P_{\varepsilon}( | \{u> \tfrac{1}{2} \} |-A)  , $$
where $\varepsilon$ is a small positive number and the penalty function is defined by
\begin{equation}\label{thepform}
 P_{\varepsilon} (s) = \left\{ \begin{array}{ll} - s/ \varepsilon \quad & \text{if} \;\; s \leq 0, \\ - \varepsilon s \qquad & \text{if} \;\; s >0 . \end{array} \right.
 \end{equation}
Penalization techniques of this type for stationary problems were analyzed in \cite{Aguilera, Tilli, Leonardi} and other works. 
The general conclusion of these treatments is that the penalized solutions converge as $\varepsilon \to 0$ to a solution of the corresponding PDE or of the original constrained minimization problem.
Moreover, it was discovered that, in order to obtain the solution of the original problem, it is not necessary to carry out the limit $\varepsilon \to 0$, since the exact minimizer is obtained for some sufficiently small but positive value of $\varepsilon$.
Although the analysis of the time-dependent problem has not yet been addressed, this observation is essential to justify the robustness of the penalty method, at least in the stationary case. \\

\noindent
{\bf Remark.}
A hint at another approach to approximating the constrained minimization can be obtained from the works developing the theory of singular perturbations \cite{Caffarelli1, Caffarelli2, Weiss}. 
The solution is realized as a limit of solutions to heat equations singularly perturbed by a nonlinear source term. 
The results also cover the time-dependent problem but are at the moment limited to free boundary conditions independent of the solution.
The application of this approach to problems of the type (\ref{nonlinheat}), where $\lambda$ is a nonlocal term depending on the solution (such as the average mean curvature of level sets), is an open problem, which we would like to address in the future.\\

In closing this section we note that, except for the paper \cite{Leonardi}, the constrained multiphase (vector-type) problem is almost unexplored. 
We provide here a formal analysis of the convergence of the BMO algorithm for the constrained multiphase problem.
Although this analysis is an important part of the present work, we defer it to \ref{app_multi} in order to make the main text less tortuous.

\section{Numerical computation}
In this section, we present a discrete algorithm for realizing constrained multiphase mean curvature flow and give the details of its numerical implementation.

\subsection{The numerical algorithm }
\label{numalg}
In treating the multiphase motions, the reformulated BMO process, stated in terms of a vector-valued heat equation as in section \ref{multinoBMO}, is approximated by use of a minimizing movement. 
A discretization in time is used to build approximate solutions by successively minimizing time-independent functionals; hence this setting conveniently allows one to include constraints via penalization. 
Each minimizer then corresponds to the solution of a vector-valued elliptic problem with Lagrange multipliers appearing on the free boundaries.

Namely, the heat equation of the BMO algorithm is solved by means of a vector-type  {\emph{discrete Morse flow}} (DMF) \cite{Rothe,Kikuchi}.
In the unconstrained case, at each step of the BMO algorithm we have to solve equation (\ref{vech}) for a time $\Delta t$.
Given regions $P_i$, $i=1, \dots ,k$, we start the DMF by constructing a function $\boldsymbol{u}_0$ such that $\boldsymbol{u}_0(x) = \boldsymbol{p}_i$ if $x \in P_i$, and set $\boldsymbol{w}_0 = \boldsymbol{u}_0$.
We choose a large positive integer $K$ which determines the time step of the flow, $h = \Delta t/K$. 
Then, to obtain the approximate interface position after time $\Delta t$, for $n=1,...,K$ we inductively minimize the following functional over $H^{1}(\Omega;{\bf{R}}^{k-1})$:
\begin{align}\label{dmf}
\mathcal{J}_{n}({{\boldsymbol{w}}}) &= \int_{\Omega}\left(\frac{|{\boldsymbol{w}}-{\boldsymbol{w}}_{n-1}|^2}{2h} + \frac{|\nabla {\boldsymbol{w}}|^2}{2} \right){dx}.
\end{align}

To evolve the interface for a time $T$, our method takes $M=T/ \Delta t$ and repeats the following for $m=1,...,M$:
\begin{enumerate}
\item Set ${{\boldsymbol{w}}}_{0}={\boldsymbol{u}}_{m-1}$. 
\item For $n=1,...,K$, compute ${\boldsymbol{w}}_{n}$ to be the minimizer of $\mathcal{J}_{n}({\boldsymbol{w}})$ over $H^{1}(\Omega;{\bf{R}}^{k-1})$.
\item Obtain $\boldsymbol{u}_m$ by thresholding:
		  	$$\boldsymbol{u}_m (x) = \boldsymbol{p}_j, \qquad \text{where} \quad \boldsymbol{p}_j \cdot \boldsymbol{w}_K(x) = \max_{i=1, \dots , k} \boldsymbol{p}_i \cdot \boldsymbol{w}_K(x) .$$

\end{enumerate}
The sequence of functions $\{ {\boldsymbol{u}}_m \}_{m=0}^M$ then gives an approximation
to the (unconstrained) multiphase motion. Figure 3 shows the basic characteristics of the method for a four-phase problem,
but see section \ref{compexp} for a clarification of the initial condition.

In treating area-constrained motions, one can see that the minimization aspect of our reformulated algorithm should allow the  inclusion of area constraints via penalization. For example, denoting the prescribed area of region $P_{i}$ by $A_{i}$, the energy functional can be modified:
\begin{equation}\label{vdmf}
	\mathcal{F}_n({\boldsymbol{w}})=\mathcal{J}_{n}({{\boldsymbol{w}}})+\frac{1}{\epsilon} \sum_{i=1}^{k} |A_{i}-meas(P_{i}^{\boldsymbol{w}})|^2.
\end{equation}
Here $\epsilon>0$ is a small penalty parameter and the areas corresponding to ${\boldsymbol{w}}$ are obtained from the sets
\begin{align}
	P_{i}^{\boldsymbol{w}} = \{ x \in \Omega ; \;\; {\boldsymbol{w}}(x)\cdot {\boldsymbol{p}}_{i} \ge  {\boldsymbol{w}}(x)\cdot {\boldsymbol{p}}_{j} \quad \forall j \}.\notag
\end{align}

Here we note that the penalty (\ref{vdmf}), which is used in our numerical computations, is slightly different from the theoretical form (\ref{thepform}). However, we prefer this form, since it is more simple and gives satisfactory results.

\begin{figure}[!ht]
\begin{center}
\includegraphics[trim=0mm 0cm 0mm 0cm,clip,height=3.75cm]{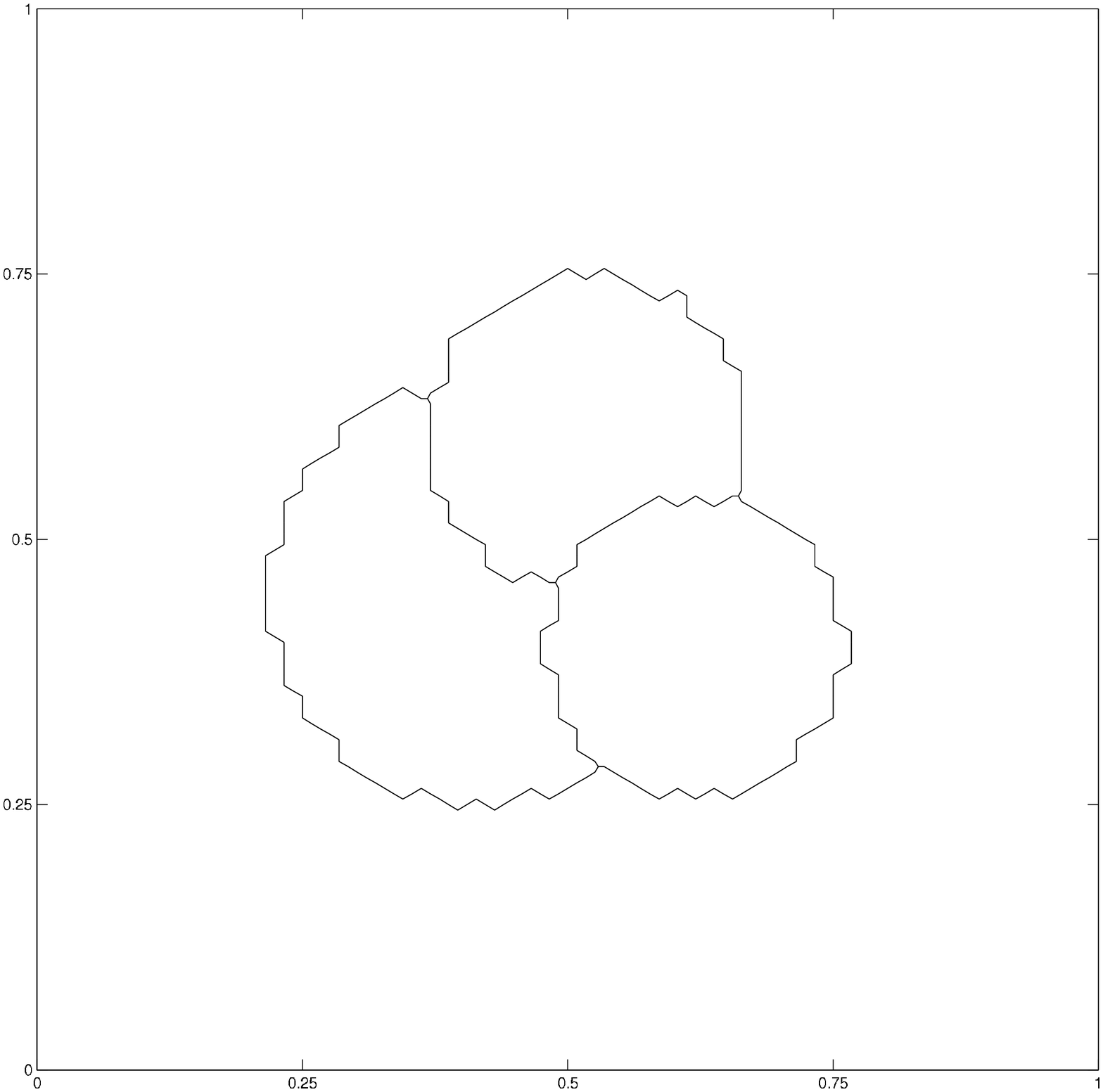}
\hspace{0.2cm}
\includegraphics[trim=0mm 0cm 0mm 0cm,clip,height=3.75cm]{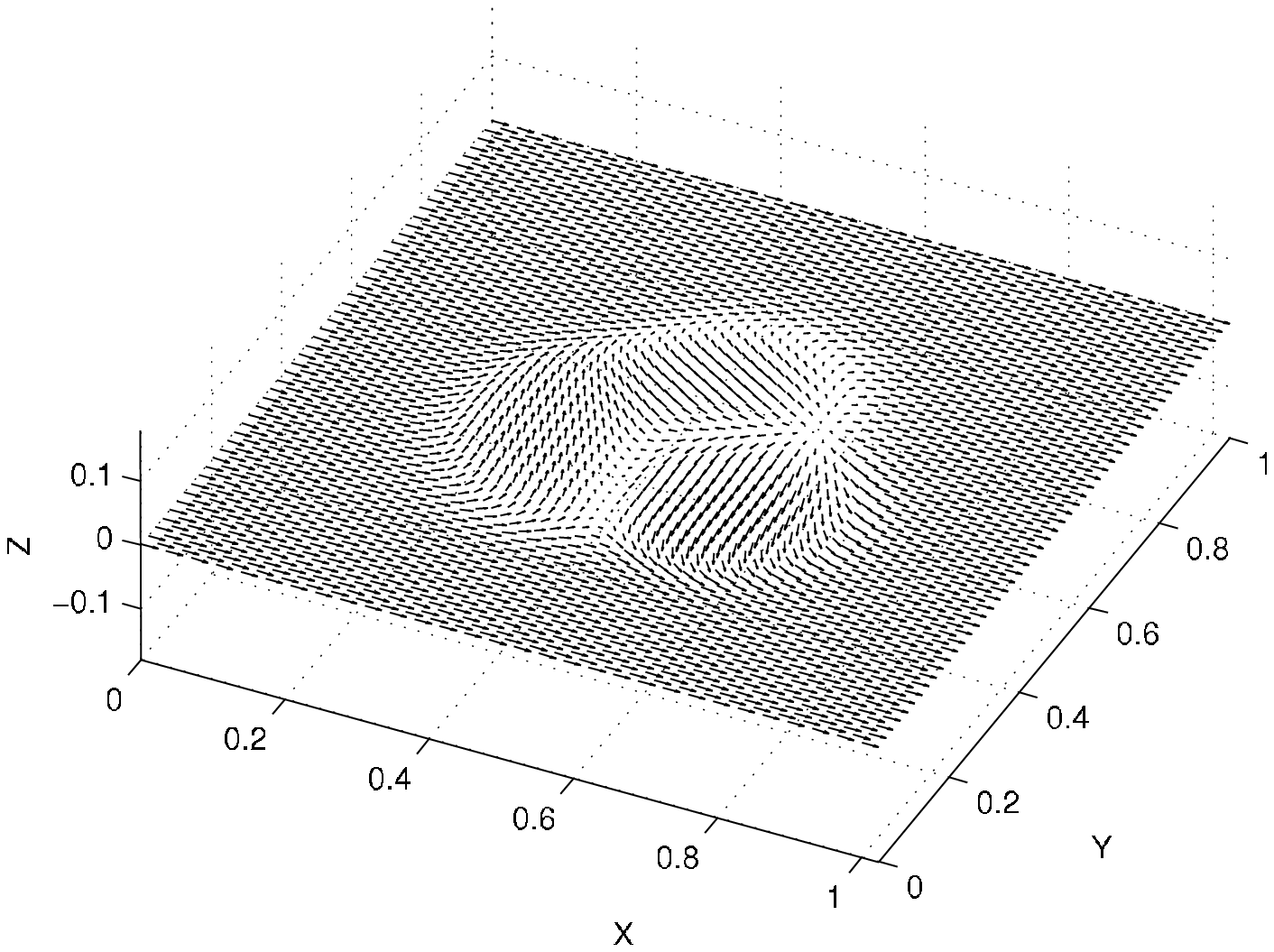}
\hspace{0.2cm}
\includegraphics[trim=0mm 0cm 0cm 0cm,clip,height=3.75cm]{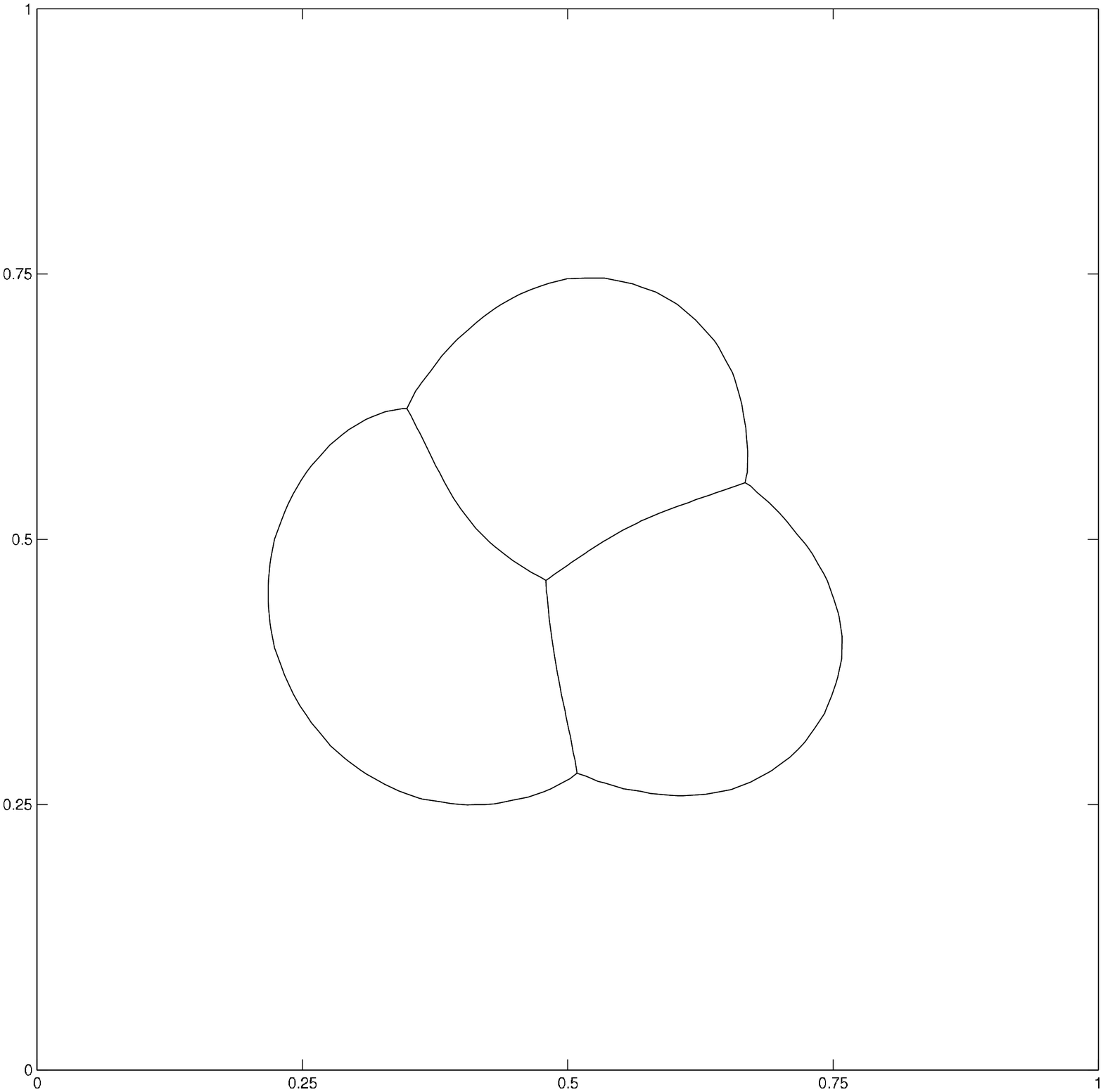}
\caption{(Left) The initial condition. (Center) An instant in time of the solution to the vector-valued heat equation. (Right) The corresponding interfaces.}
\label{fig:BMO}
\end{center}
\end{figure}

\subsection{Implementation of the method}
\label{compimp}
The numerical implementation of our method uses the finite element method to approximate the functional values (\ref{vdmf}), and minimizers are found by gradient descent.
The domain is triangulated into a finite number of elements, over which we assume that the function is continuous and a linear
interpolation between node vectors. 
The solution to the vector-type equation (\ref{nonlpdev}) is thus approximated by successive minimizations of (\ref{vdmf}) until arriving at the thresholding step. As noted in the introduction, and which is well-documented, simple thresholding by (\ref{truncstep}) is known to inhibit the motions obtained by computing with the BMO algorithm (see \cite{Esedoglu2}). We now briefly explain the troubles which may occur.

Before thresholding, the interface is a level set of the finite element solution to the heat equation, and, in the volume-preserving case, approximately satisfies the prescribed area constraint. 
However, applying the original formulation of the thresholding (\ref{truncstep}) at the nodes of the mesh would then alter the position of the interface. 
This causes two difficulties, the most significant being that, upon proceeding to the next minimizer, the interface may fail to move (thus becoming stationary). 
That is, since each evolution is obtained via a heat equation, the diffusion process must proceed long enough so that the grid resolution resolves the movement of each interface across the elements. 
On the other hand, if the diffusion proceeds too long, the approximation of the mean curvature flow looses its accuracy; hence certain configurations do not permit a suitable time step.
Additionally, due to the constraints on the area of each phase, the normal velocity of the interfaces tends to be much slower than that of the unconstrained motion, especially near the stable state. 
Therefore this issue is particularly relevant to our current problem. 
The second issue is that the enclosed areas cannot be preserved with sufficient precision after this thresholding.

We are able to overcome these issues by use of the following. 
Just prior to the thresholding step (i.e., after obtaining ${\boldsymbol{w}}_{K}$), we indicate the elements that span interfaces by $e^{*}_{j}$, record the interfacial geometry, and then threshold by (\ref{truncstep}). 
Upon the next minimization, which is the first minimization of the next BMO step, whenever we come to an indicated element, we recall the geometry of the interfaces and compute the value of the functional (\ref{vdmf}) by means of a triangulation of the element.  
In particular, the value over an indicated element is obtained by the values over a set of convex polygons, each denoted by $R_{i}$. 
These regions are determined by the element nodes and the intersection of the recorded interfaces with the element edges (see figure 4):
\begin{align}
	R_{i} =& \{x \in e^{*}_{j} ; \;\; {\boldsymbol{w}}_{K}(x)\cdot {\boldsymbol{p}}_{i} \geq  {\boldsymbol{w}}_{K}(x) \cdot {\boldsymbol{p}}_{l} \quad \forall l\ne i\}. \notag
\end{align}
For a candidate minimizer ${\boldsymbol{u}}$, the value of the discretized term in the functional (\ref{dmf}) is then computed:
\begin{equation}\label{igeo}
	\int_{e^{*}_{j}}\frac{\abs{{\boldsymbol{u}}-{\boldsymbol{u}}_{n-1}}^2}{2h}{dx} = \sum_{i=1}^{k}\int_{R_{i}\cap e^{*}_{j}}\frac{\abs{{\boldsymbol{u}}-{\boldsymbol{p}}_{i}}^2}{2h}{dx},
\end{equation}
and the contributed areas for the penalty terms
\begin{align}
	meas(e^{*}_{j}\cap P_{i}) = meas(R_{i}). \notag
\end{align}
are accumulated.
 By such an approach, we are able to realize interfacial motions whose configuration and precision of enclosed area are not altered by truncations. Moreover, this approach allows one to alleviate the restrictions on the time and space discretization of the standard BMO algorithm (see the analysis in section \ref{numtest}).
\begin{figure}[!ht]
	\begin{center}
  \includegraphics[trim=7cm 2cm 55mm 2cm,clip,height=4.2cm]{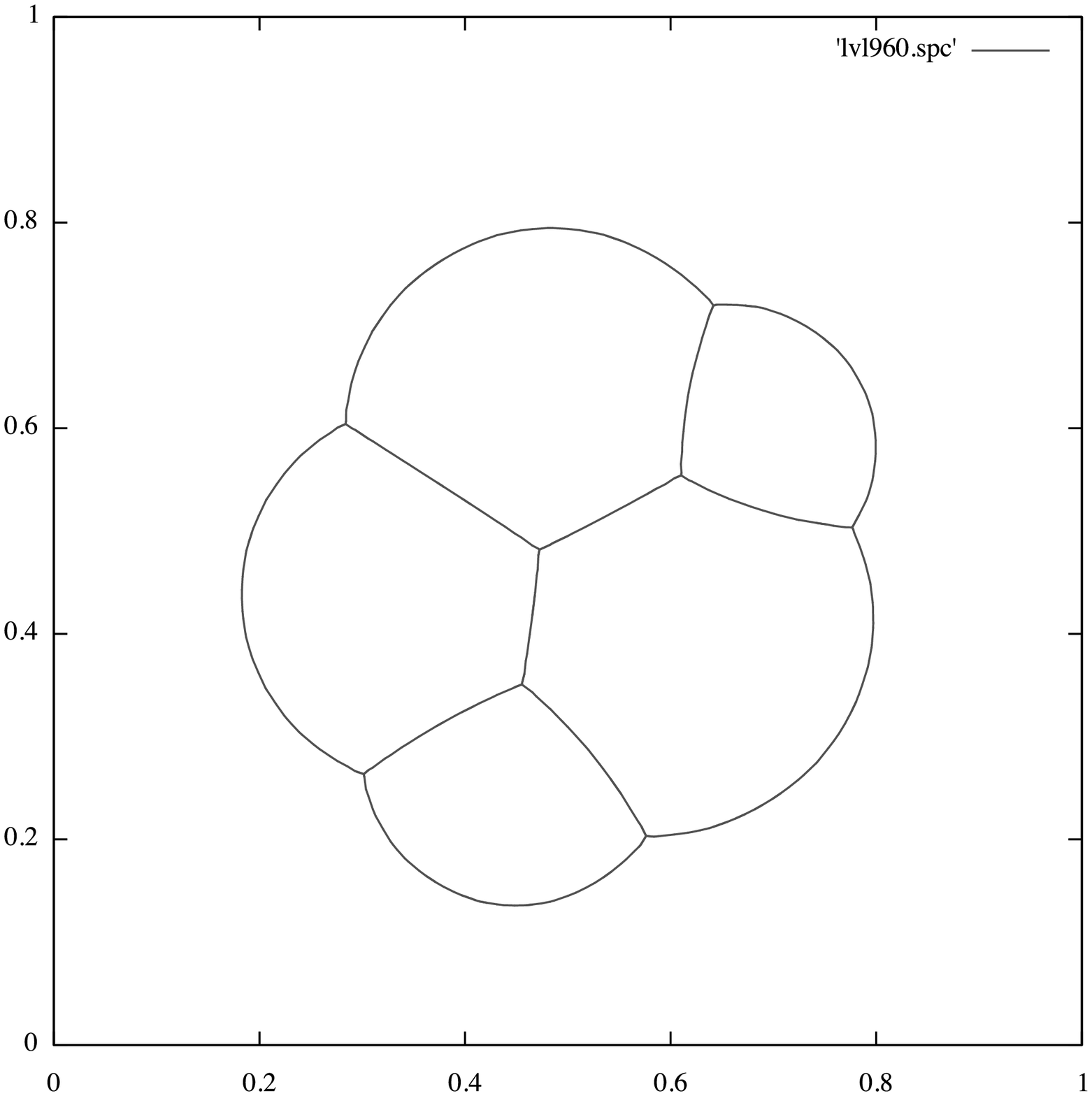}
	\includegraphics[trim=55mm 3cm 5cm 0cm,clip,height=4.2cm]{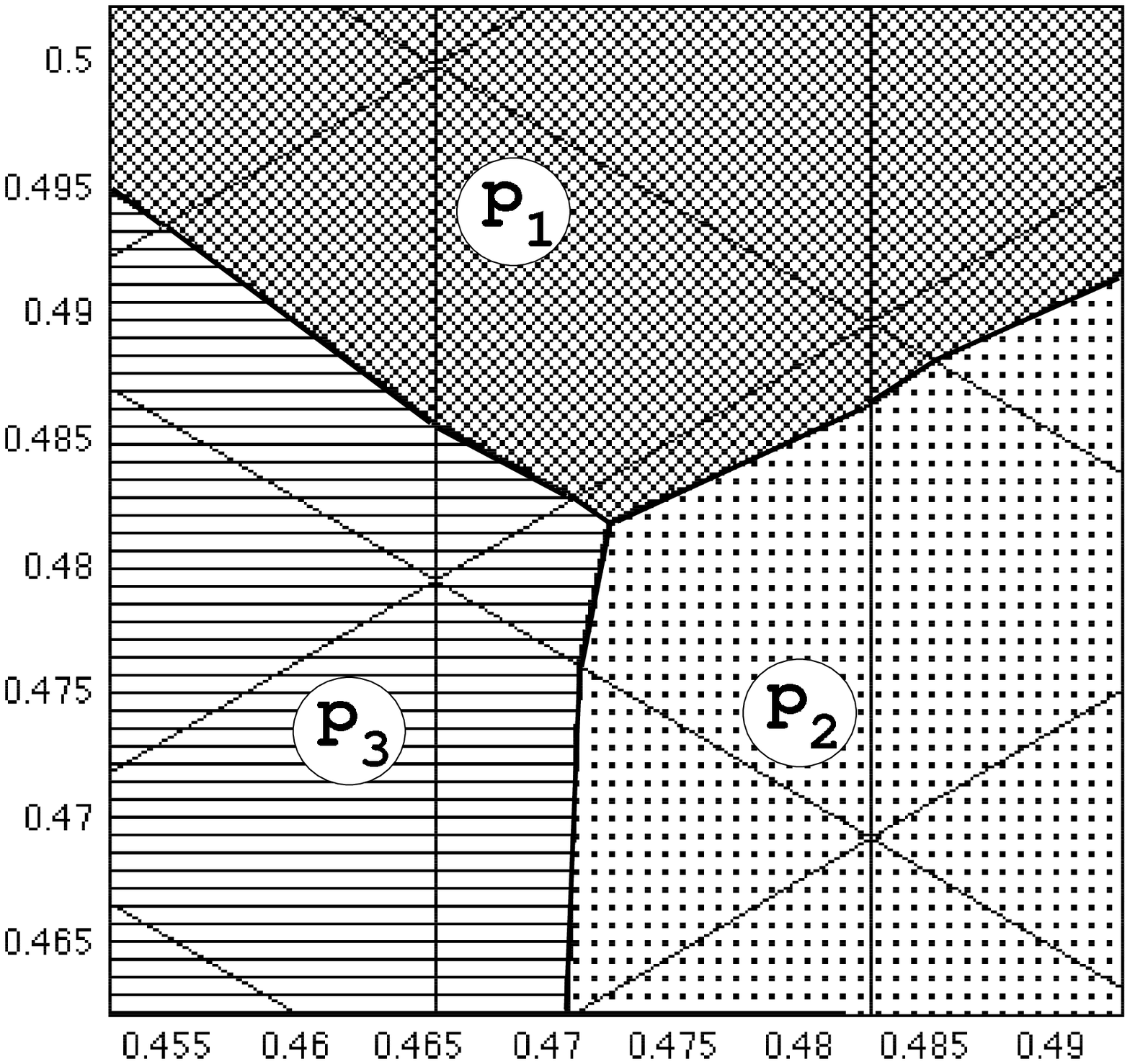} 
	\hspace{0.2cm}
	\includegraphics[scale=0.33,trim=0cm 0cm 0cm 0cm,clip,height=4.2cm]{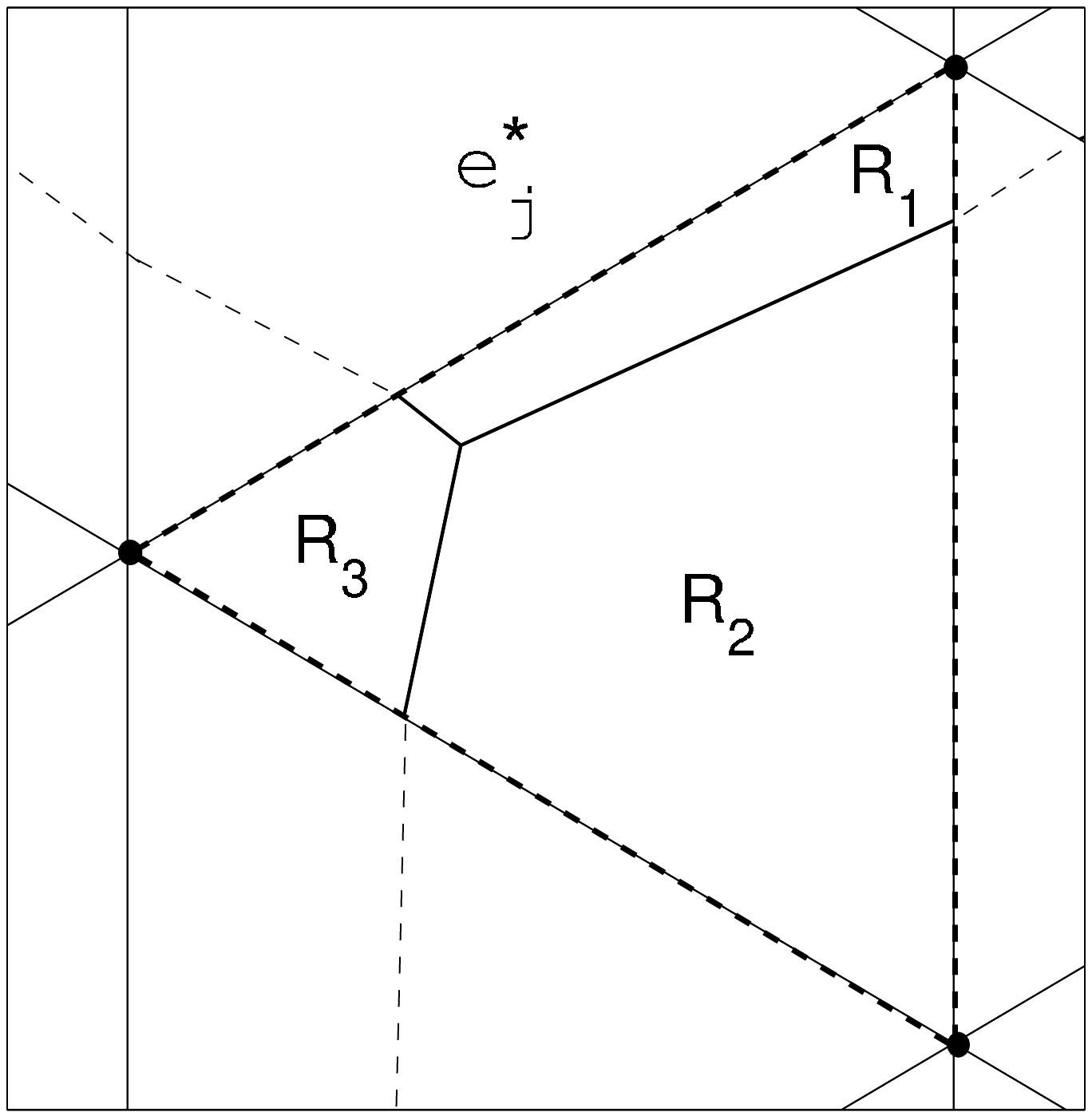}	
   	\caption{A stable solution, a junction close-up, the partition scheme for the vectors (figures use the real data).}\label{fivebub}
\end{center}
\end{figure}

\section{Numerical tests}\label{compexp}
This section presents numerical analysis and a number of numerical examples of the application of our method to area-constrained flows. 
We use ${\mathbb{P}}_{1}$ Lagrange finite elements in each of the computations. 
Thus, after assigning the appropriate vector ${\boldsymbol{p}}_{i}\in {\bf{R}}^{k-1}$ to every node of the mesh (see figure A.12 for low-dimensional visualizations of these vectors), we note the jagged shape of each initial condition. 
Although curvatures are not defined for such initial conditions, our diffusion-based algorithm is able to handle their evolution without trouble.

The physical interpretations are as follows. We configure a given number of bubbles into the shapes shown in the figures and then let them evolve. 
As we assume that there are no inertial forces, the evolution by the steepest descent (with area preservation) of the length energy (\ref{lemul}) can thus be thought of as expressing the slow movement of the bubbles. 
By examining the data, we note that the so-called symmetric Herring condition (junctions meeting at $120^{\circ}$) appears to be satisfied at junctions (see figure 4 for a close-up inspection of one such junction).

We again mention that the process described in (\ref{igeo}) is essential in computing these motions. 
Indeed, as the area-preserving interfacial evolution tends to be slower than motion by mean curvature flow, the well-known time and grid spacing restrictions of the BMO become particularly relevant \cite{Ruuth2}. 
Nevertheless, recalling the interfacial geometry after thresholding allows us to avoid such complications, and can also be used for the non-constrained BMO with the same result. 
That is, our method also works for large ratios of grid and time step sizes, for which the original BMO becomes stationary. 
Moreover, formal analysis and numerical tests show that this approach does not alter the characteristics of the target motion.  

It is also worth stating that another technique for handling the restriction of the BMO algorithm on non-adaptive meshes is through the use of signed distance functions. 
This method was developed in \cite{Esedoglu2}, where it was shown that it gives satisfactory results and we would like to extend the method used there to the multiphase constrained case.

\subsection{Convergence analysis}

\subsubsection{Analysis of error}
\label{numtest}
This section examines the behavior of our method in comparison to the standard BMO. We refer to the standard method as BMO, and to our own algorithm (which utilizes the thresholding from section \ref{compimp}) as BMO$^{*}$.

We examine the error of our method when applied to a simple test problem. 
By symmetry, a circle of initial radius $r_{0}$ which is evolving by mean curvature flow remains a circle whose radius $r(t)$ satisfies the following ordinary differential equation:
$$ r' = -\frac{1}{r} \hspace{20pt}\text{(solution: $r(t)=(r_{0}^2-2t)^{1/2}$)}. $$
With an initial condition obtained from the target radius $r_{0}=0.35$, we vary grid and BMO time step sizes, and compute the solution by use of the BMO and BMO$^{*}$ algorithms. 
The output of the program is a list of interface nodes $\{(P^x_i,P^y_i)\}_i$.
Circles were fitted to these points at each time level by minimization of the functional
$$ \sum_i \Big( (P^x_i-C^x)^2 + (P^y_i-C^y)^2-r^2 \Big)^2 $$
with respect to the centre coordinates $(C^x,C^y)$ and radii $r$.
We measure the error of the method by the time-average of the absolute difference between the radius of the fitted circle and the exact radius:
$$ \frac{1}{L} \sum_{l=1}^L | r_{\text{fit}}(t_l)-r(t_l)| . $$
Here $L$ denotes the number of time steps until the radius is zero.
The error table for the standard BMO algorithm is as follows:

{\footnotesize
{\begin{equation}
\begin{array}{|c|c|c|c|c|c|c|c|c|}
\hline
\text{res(space$\backslash$time)}&2 & 4 & 8 & 16 & 32 & 64 & 128 & 256\\ \hline
5 \times 5&0.0483 & 0.0231 & - & - & - & - & - & - \\
10 \times 10&0.0572 & 0.0176 & - & - & - & - & - & -\\
20 \times 20&0.0047 & 0.0018 & 0.0166 & 0.0070 & - & - & - & - \\
40 \times 40&0.0044 & 0.0039 & 0.0034 & 0.0036 & - & - & - & -\\
80 \times 80&0.0056 & 0.0032 & 0.0027 & 0.0101 & 0.0042 & 0.0223 & - & -\\
160 \times 160&0.0060 & 0.0031 & 0.0028 & 0.0073 & 0.0107 & 0.0039 & 0.0040 & 0.0397\\
\hline
\end{array}\notag
\end{equation}}
}

We confirm that when the time step decreases to a certain size (relative to the gird), the BMO method halts (this is indicated by the symbol ``--").
The critical ratio of space to time step size is approximately $20$.

On the other hand, by computing the same evolutions with the BMO* algorithm, 
we find that it is able to deal with a much wider range of parameters (critical ratio around $400$).
The initial condition is an approximation to a circle of radius $r_0=0.35$ obtained from the area-preserving BMO$^*$ method with penalty parameter $\epsilon = 10^{-6}$.
The error table is as follows:
{\footnotesize
{\begin{equation}
\begin{array}{|c|c|c|c|c|c|c|c|c|}
\hline
\text{res(space$\backslash$time)}&2 & 4 & 8 & 16 & 32 & 64 & 128 & 256\\ \hline
5 \times 5&0.0276 & 0.0228 & 0.0258 & 0.0381 & 0.0440 & 0.0368 & - & -\\
10 \times 10&0.0121 & 0.0146 & 0.0076 & 0.0165 & 0.0201 & 0.0103 & 0.1210 & -\\
20 \times 20&0.0044 & 0.0038 & 0.0033 & 0.0043 & 0.0081 & 0.0139 & 0.0080 & 0.0694\\
40 \times 40&0.0045 & 0.0035 & 0.0024 & 0.0024 & 0.0028 & 0.0007 & 0.0056 & 0.0096\\
80 \times 80&0.0053 & 0.0032 & 0.0025 & 0.0042 & 0.0070 & 0.0073 & 0.0056 & 0.0016\\
160 \times 160&0.0059 & 0.0031 & 0.0026 & 0.0049 & 0.0085 & 0.0097 & 0.0097 & 0.0088\\
\hline
\end{array}\notag
\end{equation}}
}

Overall, we see that both algorithms approximate the solution, and that the additional partitioning of the BMO$^{*}$ algorithm is beneficial. 
For finer meshes, both algorithms show a tendency of stagnating errors when the space mesh is refined.
This fact could be attributed to the properties of the DMF scheme.
However, we note that the point of this partitioning is not to improve the error per se, but to relax the time and grid restrictions inherent to the standard BMO.

\subsubsection{Analysis of penalty parameter}
Here we perform an error analysis for the two-phase area-preserving case. 
Since the velocity of the area-preserving motions tends to be much slower than in the non-constrained case, the use of the BMO$^{*}$ should be preferred over the standard BMO.

We take two non-intersecting circles of radii $r_{a}=0.1996$ and $r_{b}=0.1384$ and identify them as the same phase. 
Then the area preservation condition implies that the radius of the larger circle will grow as the smaller circle's radius shrinks.
The evolution of the radii follows the equations
\begin{align}
r'_{1} &= \frac{-1}{r_{1}} + \frac{2}{r_{1}+r_{2}} \notag\\
r'_{2} &= \frac{-1}{r_{2}} + \frac{2}{r_{1}+r_{2}} \notag.
\end{align}
We use a numerical method to compute a precise approximation to the solution of the above system and compare the result to that obtained by our penalty method.
In the BMO* method, the domain $[0,1]\times[0,1]$ is triangulated into 14536 elements, $\Delta t = 2.5\times10^{-4}$, and $K=30$. 
The results are ploted in figure 5. 

Weak penalties with $\epsilon = 10^0, 10^{-1}, 10^{-2}$ give almost the same result, and the larger circle shrinks slightly although it should grow.
For increasing strength of the penalty the results approach the correct solution with the larger circle growing as the smaller shrinks.
Finally, penalties $\epsilon = 10^{-5},10^{-6},10^{-7}$ give almost identical results, which are only slightly different from the exact solution.
This finding agrees with the theoretical prediction mentioned at the end of section \ref{mmac}; namely that we should obtain the exact
solution for sufficiently large penalties.
The solutions for large penalties indeed do not change, and their slight deviation from the exact solution is caused by the discretization.
In conclusion we can say that the penalty method behaves well, since for each grid resolution there exists a range of penalty coefficients for which the solution is independent of the penalty strength and appropriately approximates the exact solution.

\noindent
\begin{figure}[!ht]
\begin{center}
\includegraphics[trim=0mm 0cm 0mm 0cm,clip,height=7.0cm]{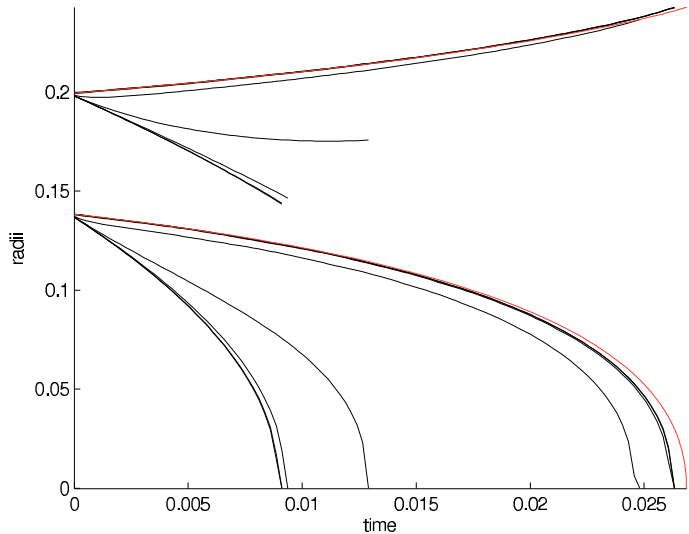}
\caption{Evolution of the radii for penalty parameters $\epsilon$ varying from $10^0$ to $10^{-7}$ (black curves, ordered from left to right). The red curve shows the exact solution.}
\label{fig:radii}
\end{center}
\end{figure}

\subsubsection{Analysis of the multiphase area-preserving algorithm}

In order to test the performance of our multiphase area-preserving algorithm, we compute the stationary solution corresponding to two 2-dimen\-sio\-nal soap bubbles attached to a wall. 
It can be rigorously shown that the steady shape is composed of three circular arcs meeting with $120^{\circ}$ angles at the triple junction and with $90^{\circ}$ angles at the walls. 
Moreover, the radii of the arcs satisfy the well-known condition 
$$ \frac{1}{r_1} - \frac{1}{r_2} = \frac{1}{r_{12}} ,$$
where $r_1$ and $r_2$ are the radii of the bubbles and $r_{12}$ is the radius of their common arc.
When the initial volumes of the bubbles are given, it is possible to compute the above radii analytically.

This analytical solution is compared to the numerical solution in figure 6.
The numerical solution was obtained by running the area-preserving three-phase method for sufficiently long time, until the interfaces stopped moving.
The domain $[0,1]\times[0,1]$ was triangulated into 14536 elements, $\Delta t = 2.5\times10^{-4}$, $K=30$, $\epsilon = 10^{-6}$, 
and the fitting method from section \ref{numtest} was used to obtain the circle radii.
Although the area differs slightly from the prescribed value, the configuration of the numerical solution, including the triple junction, agrees well with the analytic solution relative to the resolution of the grid.  
\noindent
\begin{figure}[!ht]
\begin{center}
\includegraphics[trim=0mm 0cm 0mm 0cm,clip,height=5cm]{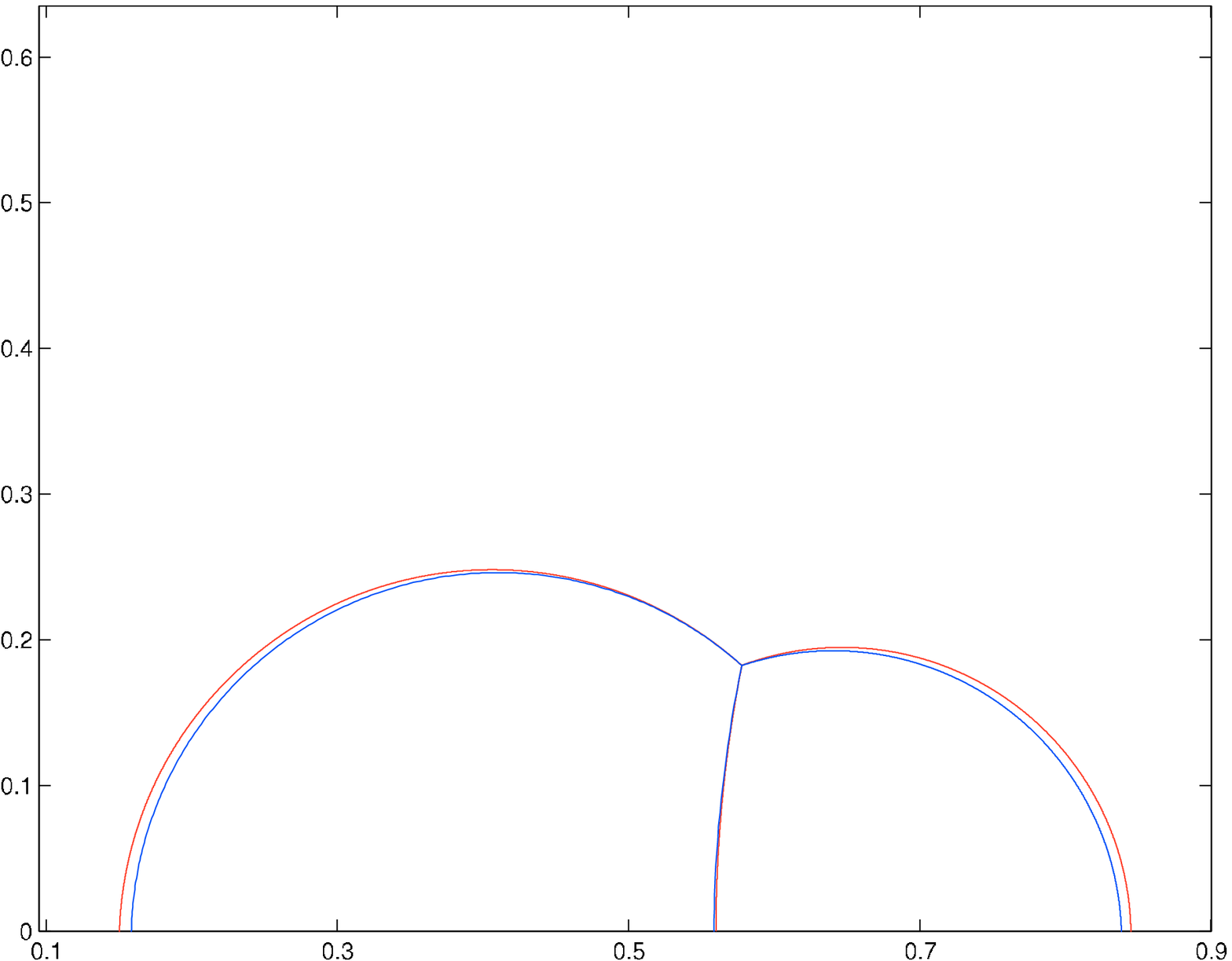}
\hspace{0.3cm}
\includegraphics[trim=0mm 0cm 0mm 0cm,clip,height=5cm]{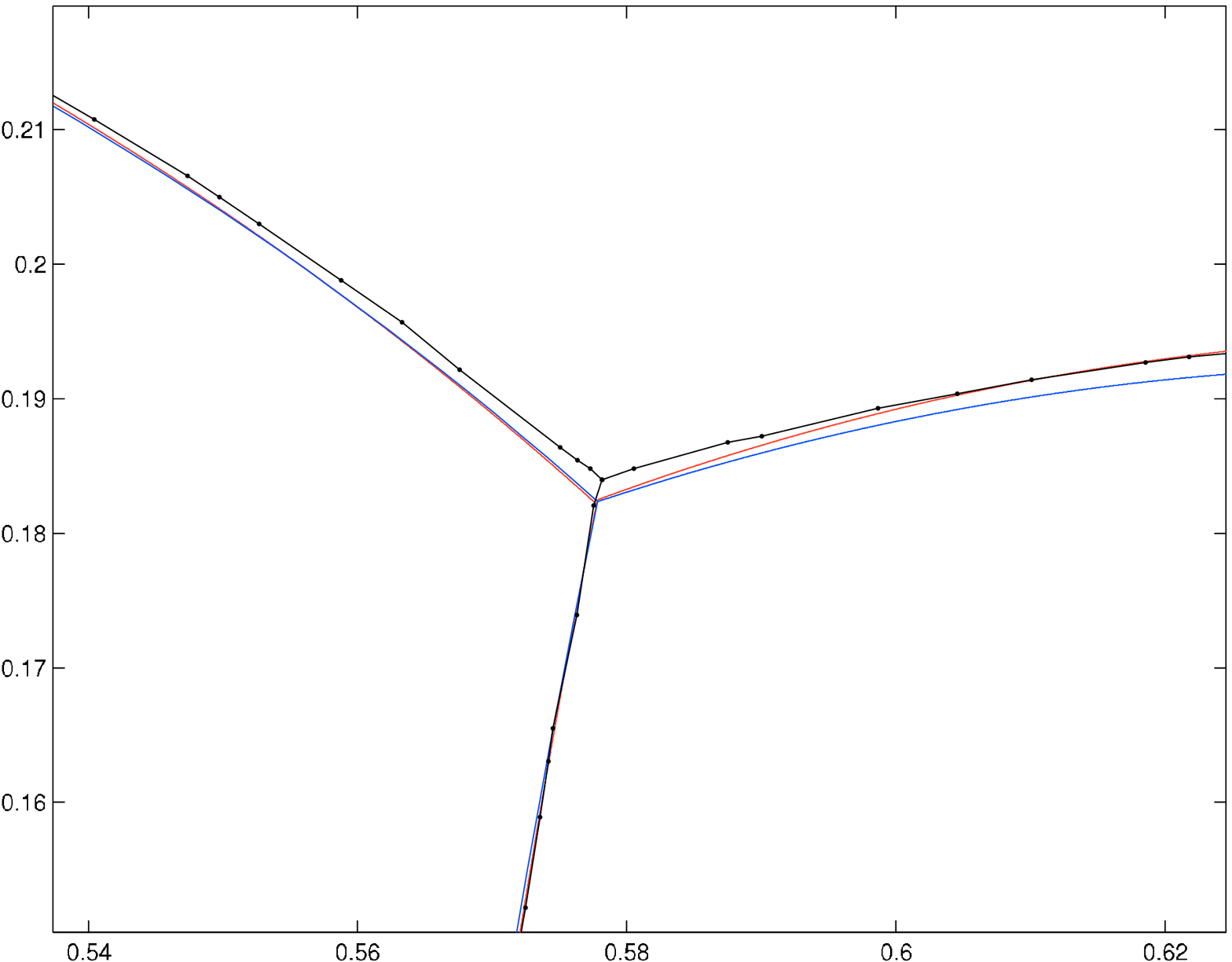}
\caption{The stationary configuration for two bubbles attached to a wall: (Left) Exact solution (blue) and numerical solution (red) -- circles obtained by least squares fitting are plotted. (Right) A close-up view of the triple junction showing the exact solution (blue), fitted numerical solution (red) and original numerical solution (black with dots).}
\label{fig:multistat}
\end{center}
\end{figure}

\subsection{Multiphase flow}
Here we present two examples of multiphase mean curvature flows.

\subsubsection{A triple bubble}\label{stdparam}
We begin by examining the motion of four phases, where the three bubbles initially touch each other (figure 7). 
The area of each phase is different, and a corresponding triple bubble is obtained as the stable configuration for large times.

Here the domain $\Omega=[0,1]\times[0,1]$ is triangulated into approximately 5600 elements, $\Delta t = 5\times10^{-4}$, and $K=10$ . 
The mesh is nearly uniform so that most elements have area approximately equal to $1.8 \times 10^{-4}$. 
A penalty of the form shown in (\ref{vdmf}) is added for each phase and its parameter is $\epsilon = 10^{-6}$. 
The prescribed volumes were maintained to within an absolute error of $10^{-4}$, even during the dynamic portion of the movement.

\noindent
\begin{figure}[!ht]
\begin{center}
\includegraphics[trim=0mm 0cm 0mm 0cm,clip,height=4.0cm]{tripvpwlinit.eps}
\hspace{0.2cm}
\includegraphics[trim=0mm 0cm 0mm 0cm,clip,height=4.0cm]{tripvpwl.eps}
\hspace{0.2cm}
\includegraphics[trim=0mm 0cm 0cm 0cm,clip,height=4.0cm]{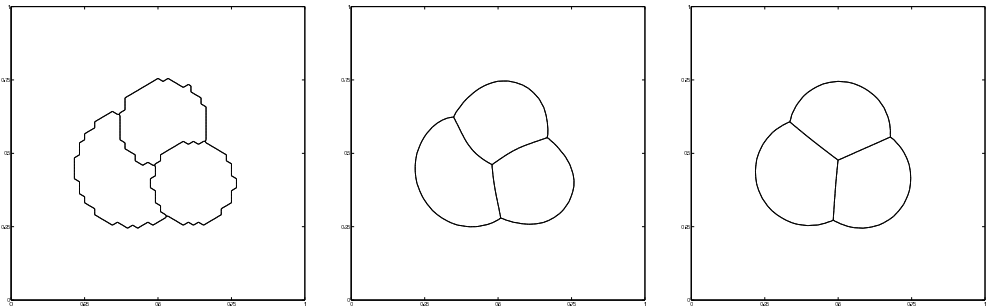}
\caption{(Left) The initial condition. (Center) Evolution after a short time. (Right) The stable configuration.}
\label{fig:trip}
\end{center}
\end{figure}

\subsubsection{A nine phase flow}
In the present computation, nine phases are positioned throughout the domain.  
After assigning the appropriate vector to each node location and detecting the interfaces, we have the initial interfaces as shown in figure 8. 
The bubbles cling to the lower boundary and, eventually, the top most bubble slides itself in between two others. 
This shows that the method can naturally handle topological changes.

The domain $\Omega = [0,1]\times[0,1]$ is triangulated into approximately 9600 elements, $\Delta t=3 \times 10^{-4}$, and $K=30$. 
The mesh is nearly uniform so that most elements have area approximately equal to $10^{-4}$. 
A penalty of the form shown in (\ref{vdmf}) is again added for each phase and its parameter is $\epsilon =10^{-6}$. 
The absolute errors in the areas were able to be maintained within $10^{-4}$.
\begin{figure}[htbp]
	\begin{center}
	\includegraphics[height=4.0cm]{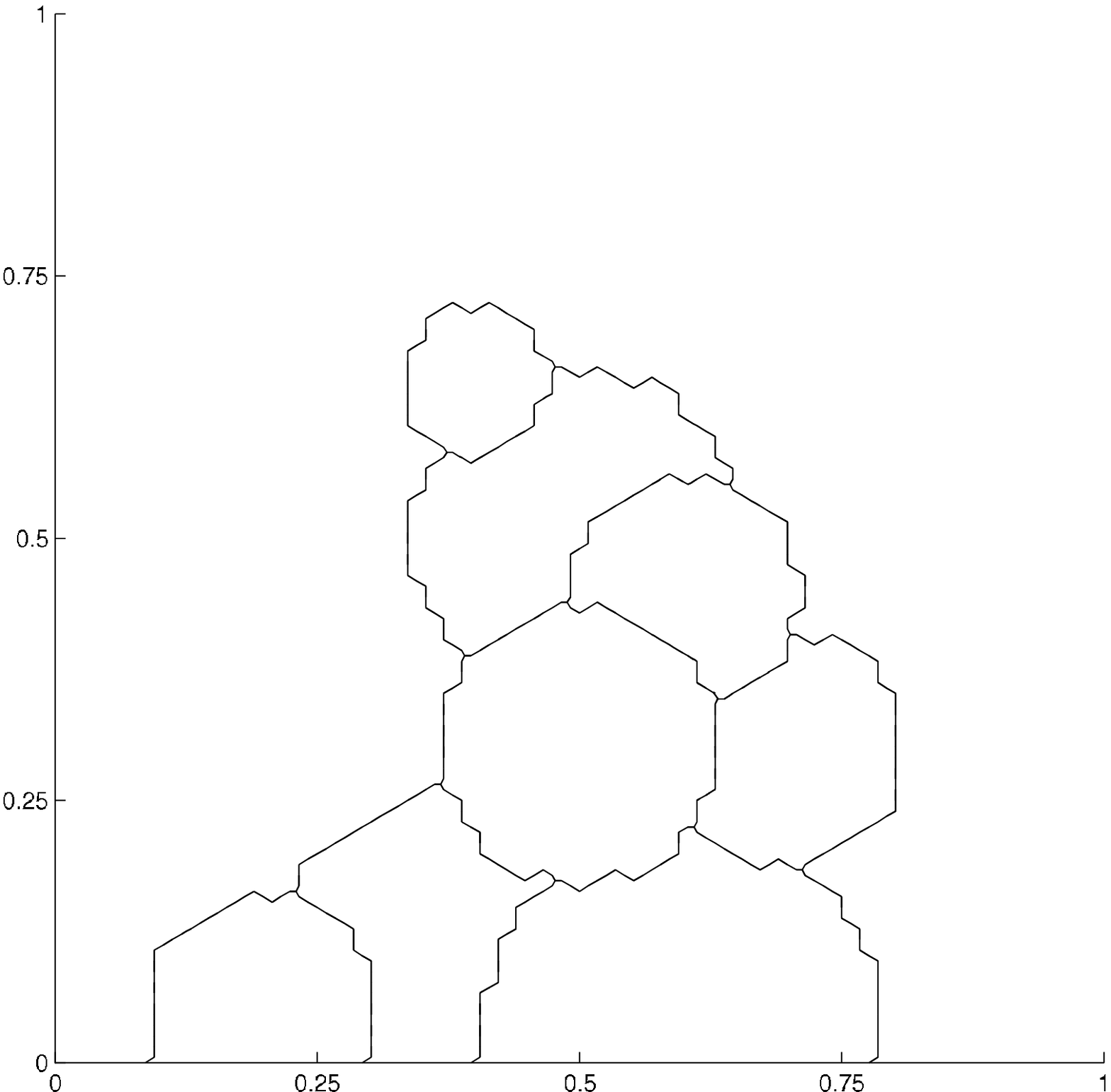}
	\hspace{0.2cm}
	\includegraphics[height=4.0cm]{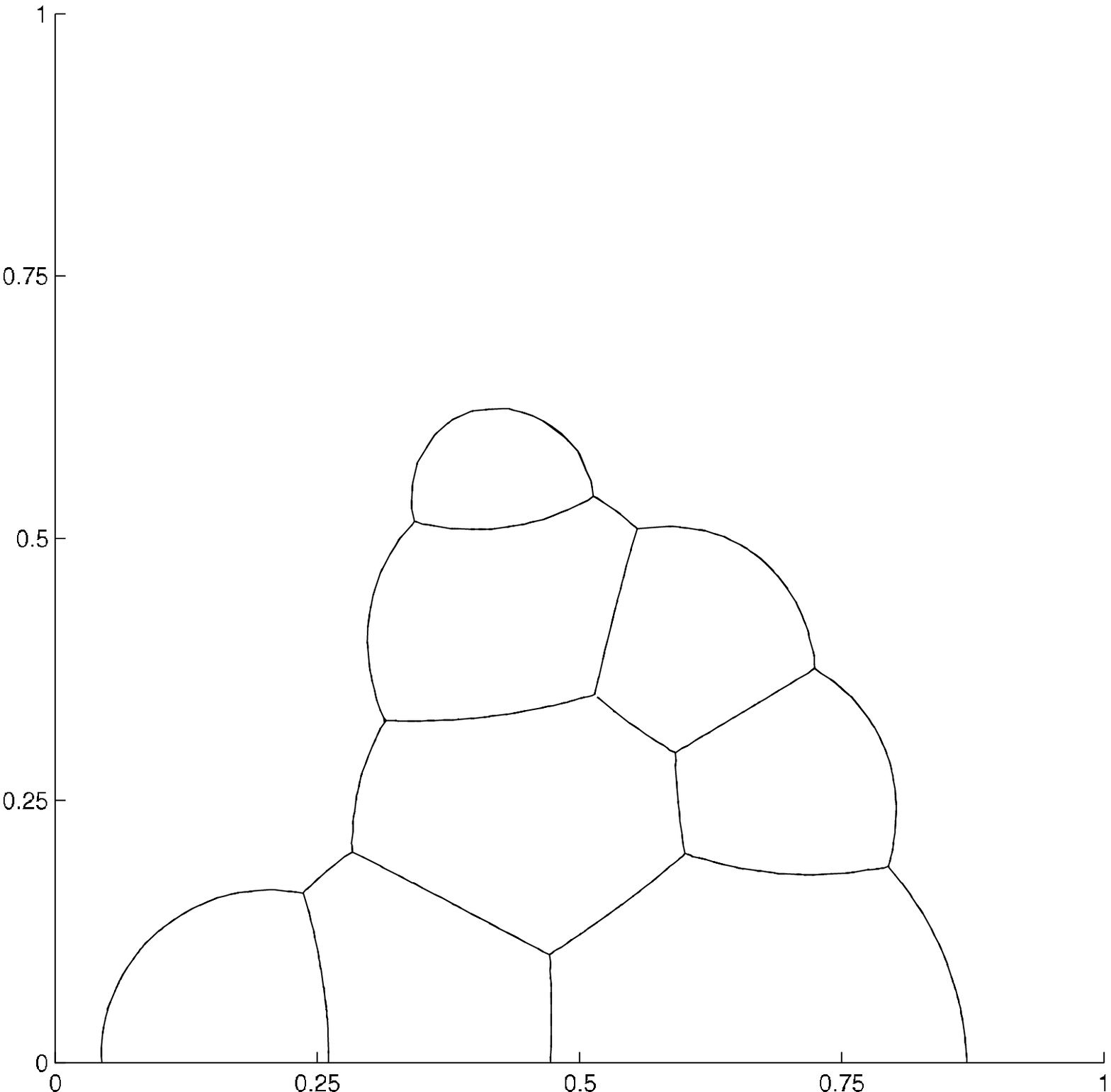} 
	\hspace{0.2cm}
	\includegraphics[height=4.0cm]{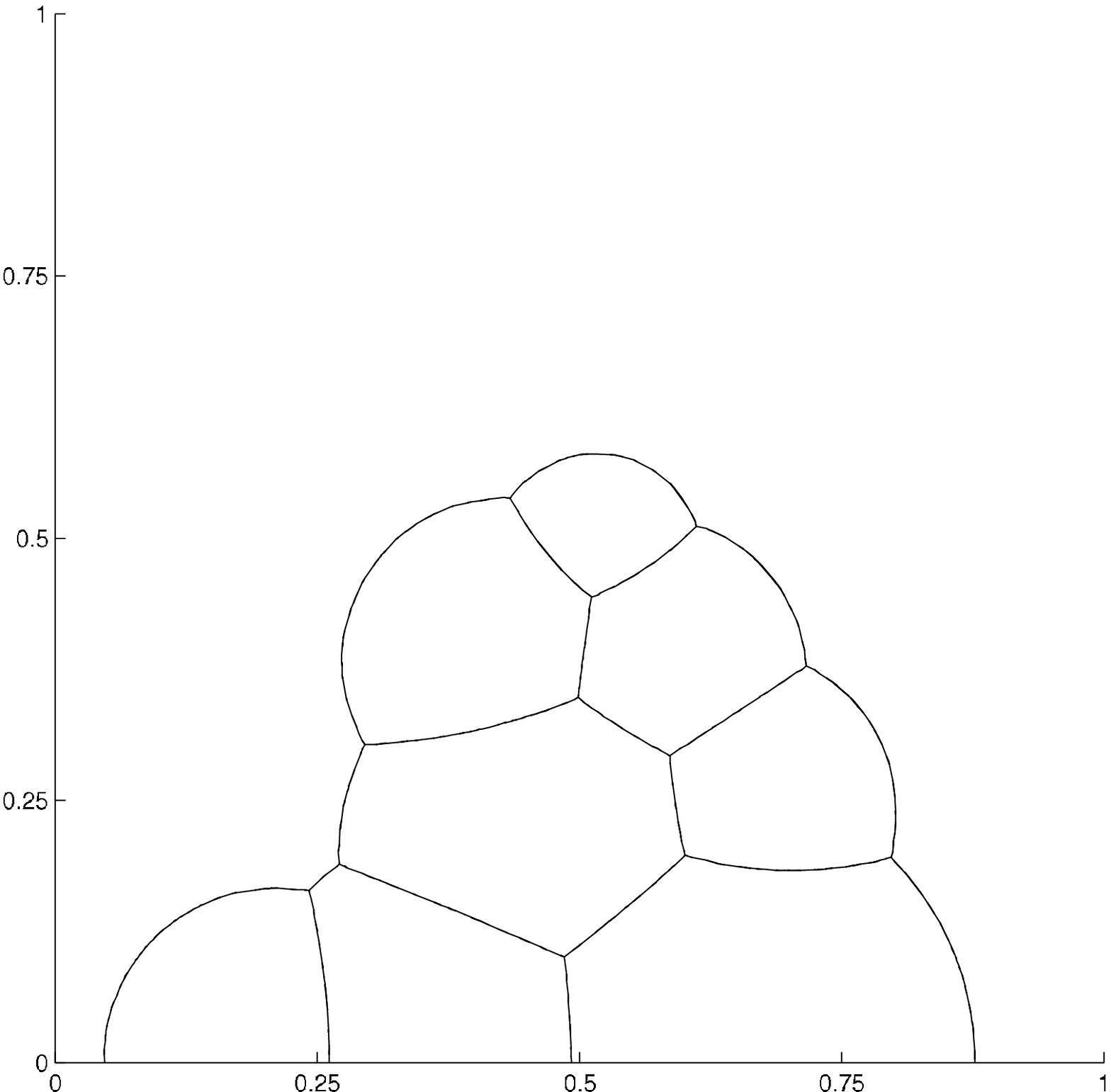}	
	 	\caption{The initial condition, evolution after $500\Delta t$, and the stable configuration.}\label{ninebub}
	 	\end{center}
\end{figure}

\subsection{Interaction of interfaces}
The proposed method can also handle interactions of interfaces, such as merging and attaching.

\subsubsection{Two-phase coalescence}
We first consider a computation involving two phases. 
One phase is initially separated into two distinct parts and is configured in such a way that, in the course of the interface evolution, the two parts eventually come into contact with each other. 
The parameters are as in section \ref{stdparam}.

When the two parts touch, a topological change occurs, and the evolution continues until reaching a stable configuration (a circle); see figure 9. 
We note that, due to the diffusion process of our method, the initially separate phases attract each other slightly. 
Nevertheless, this attraction decreases extremely fast with the distance of the phases.
However, for this reason, capturing the precise behavior at the moment of the topological change is difficult. 
Of course, this unwanted attraction can be reduced by refining the time step and grid resolution.
The prescribed volumes were maintained to within an absolute error of $10^{-4}$.
\begin{figure}[!ht]
	\begin{center}
	\includegraphics[width=4cm]{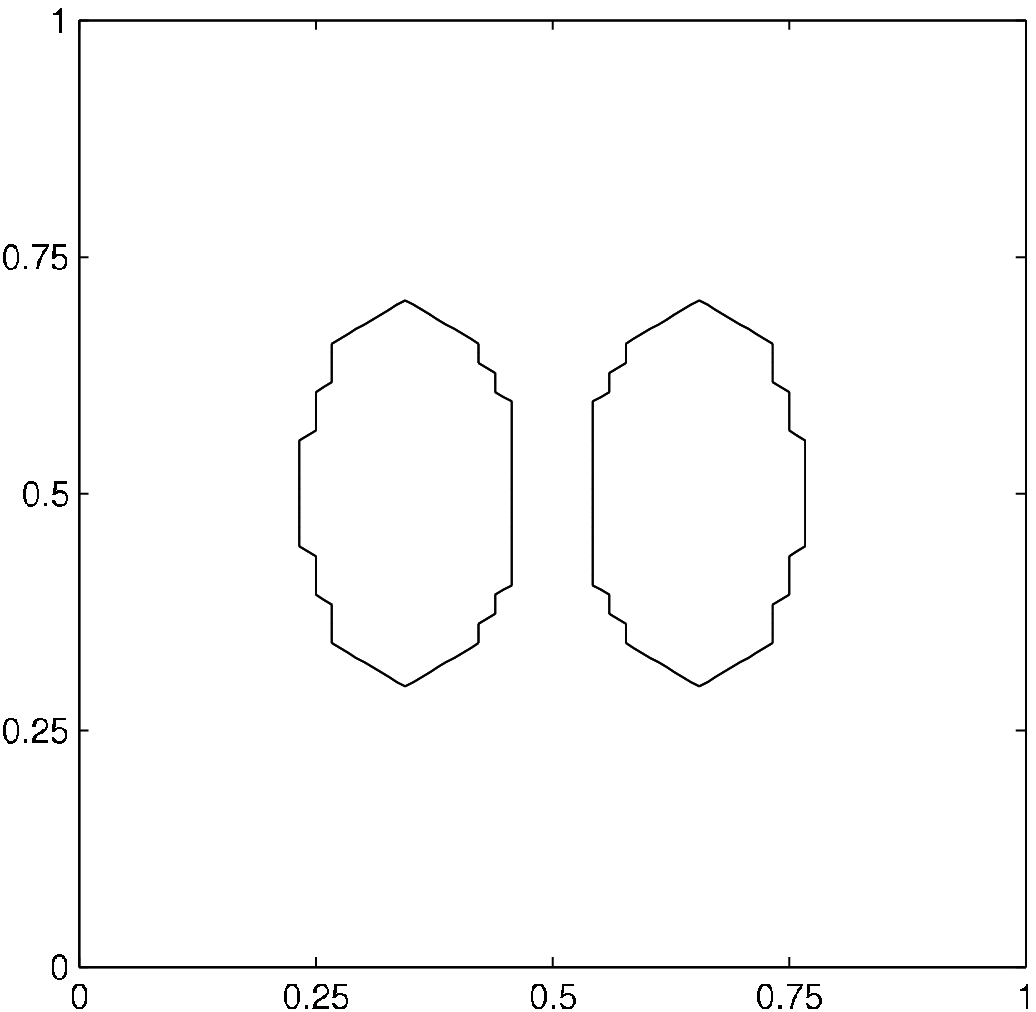} \hspace{1cm}
	\includegraphics[width=4cm]{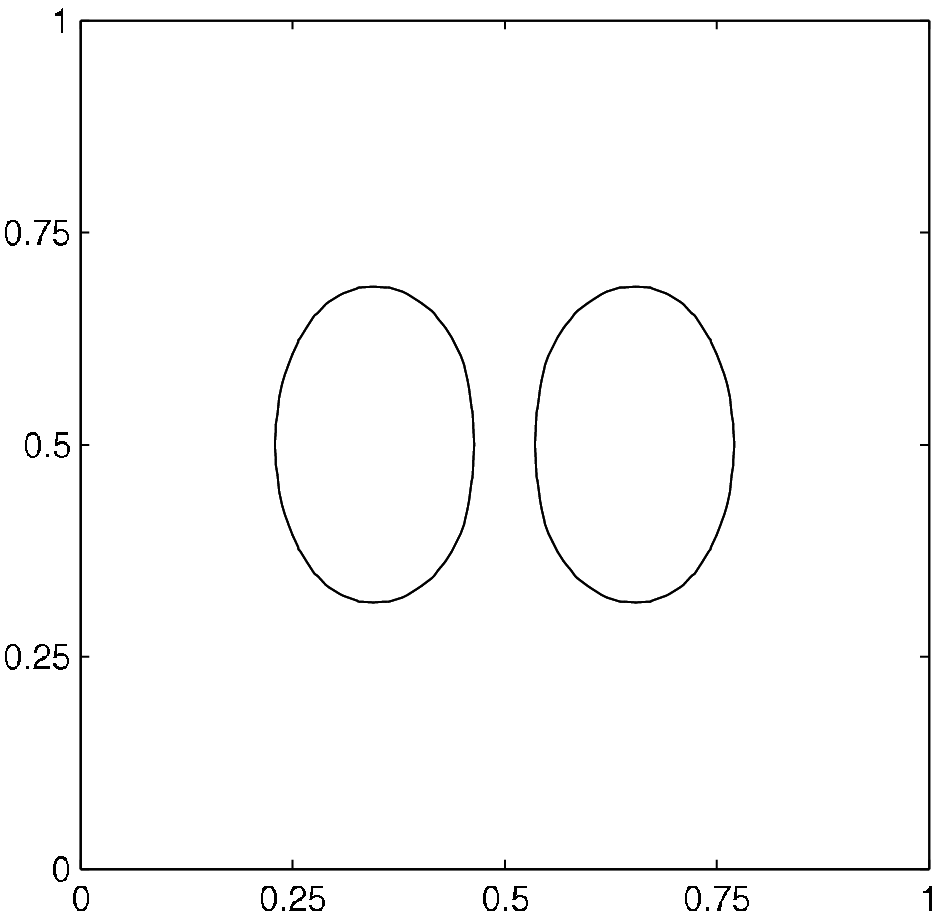}	\\
	\includegraphics[width=4cm]{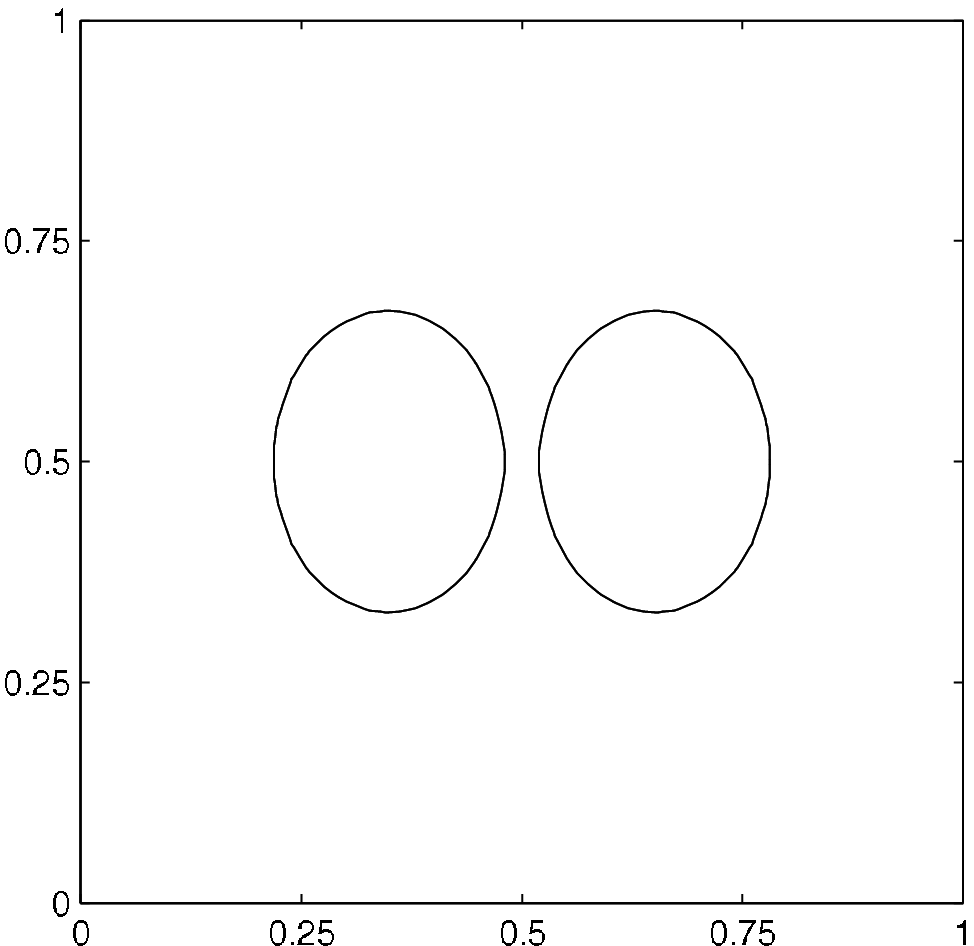} \hspace{1cm}
	\includegraphics[width=4cm]{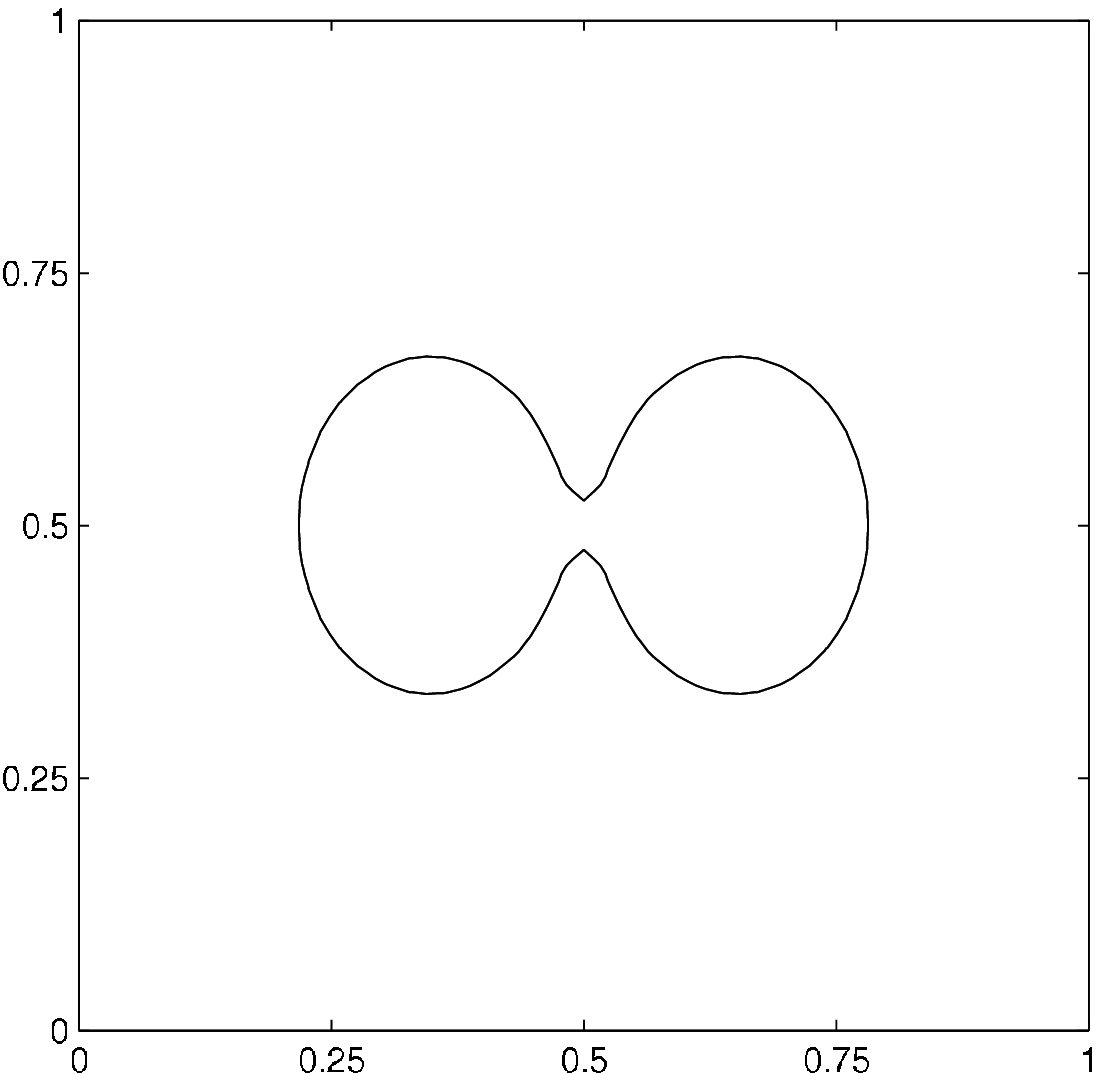} \\
	\hspace{0.12cm}
	\includegraphics[width=3.8cm]{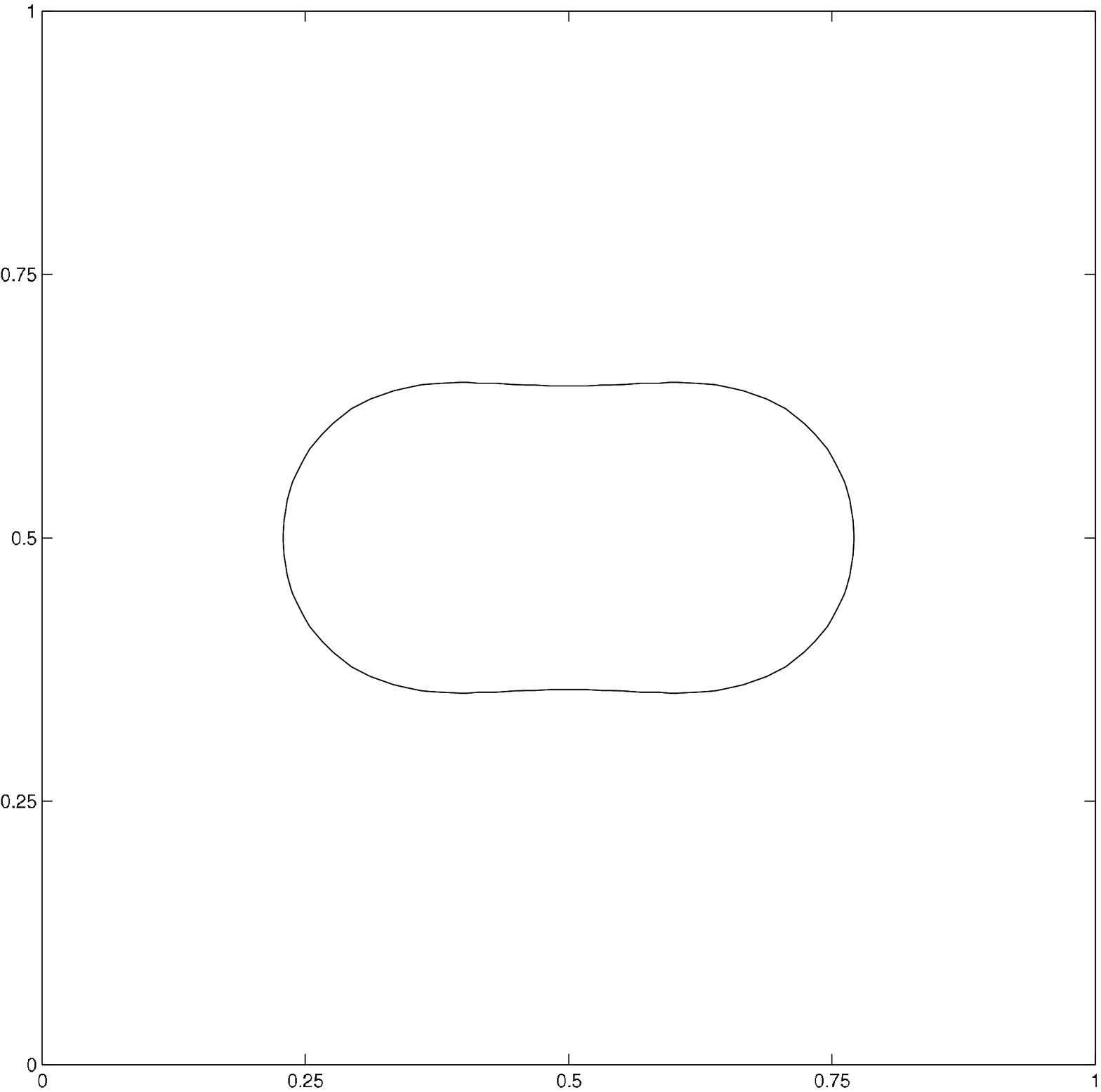} \hspace{1.03cm}
	\includegraphics[width=4cm]{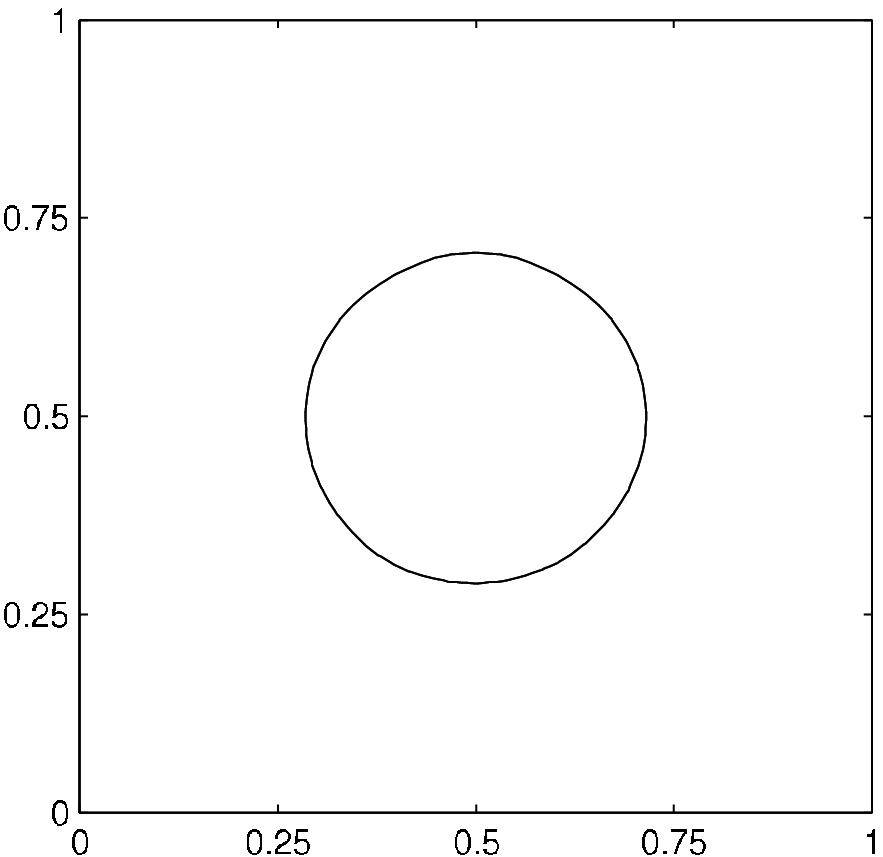}	
	 	\caption{Merging of two regions corresponding to the same phase.}\label{oval2p1}
\end{center}
\end{figure}

\subsubsection{Three phase coalescence}
In this computation, we place three phases throughout the domain. 
Two are configured to be initially separate, but so that they eventually touch. 
As the areas of the phases are different, the final steady state solution shows the shape of a non-symmetric double-bubble; see figure 10. 
The setting of parameters is identical to that of section \ref{stdparam}.
Absolute errors in the areas are approximately $10^{-4}$. 
\begin{figure}[!ht]
	\begin{center}
	\includegraphics[width=4cm]{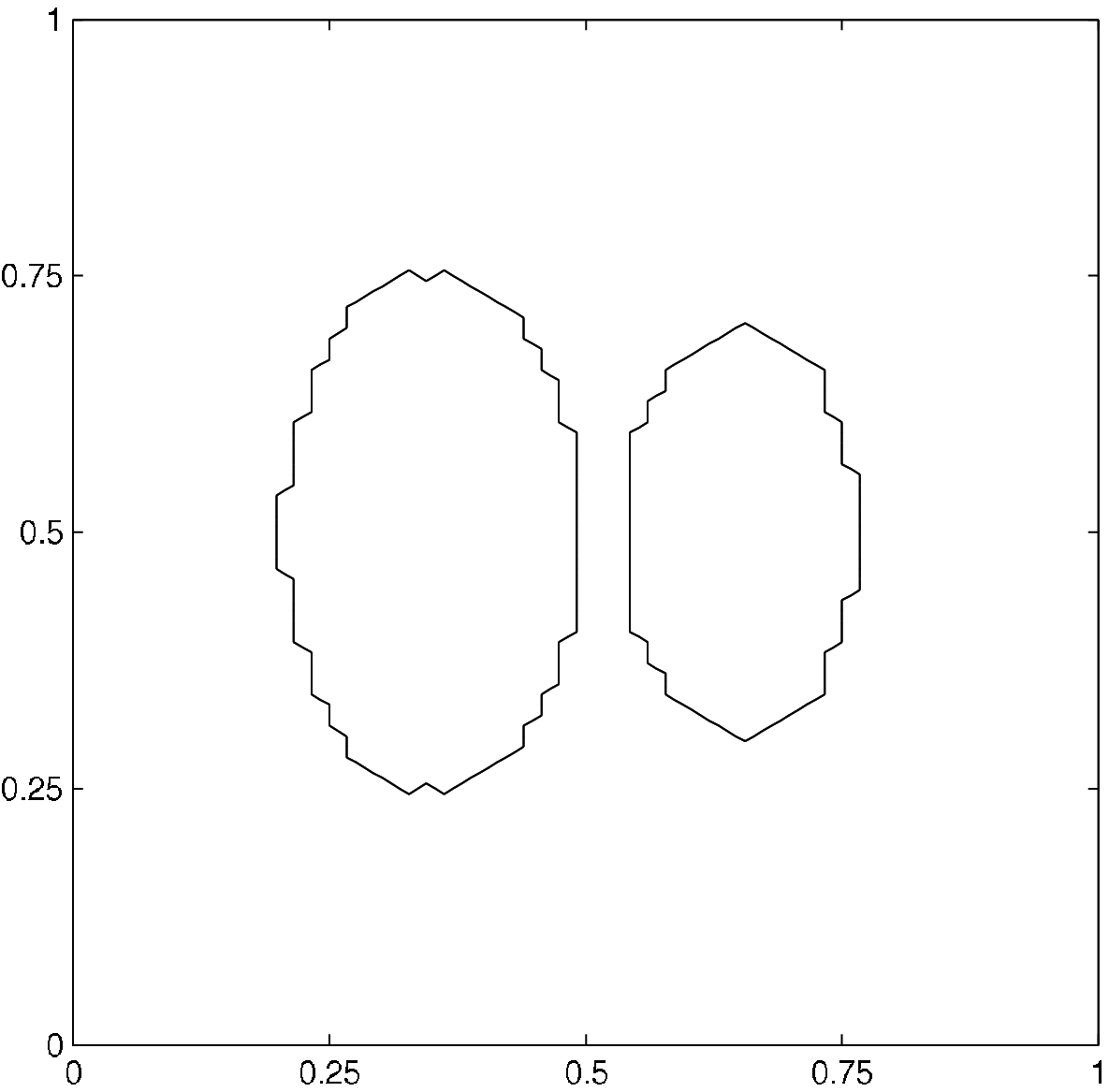} \hspace{1cm}
	\includegraphics[width=4cm]{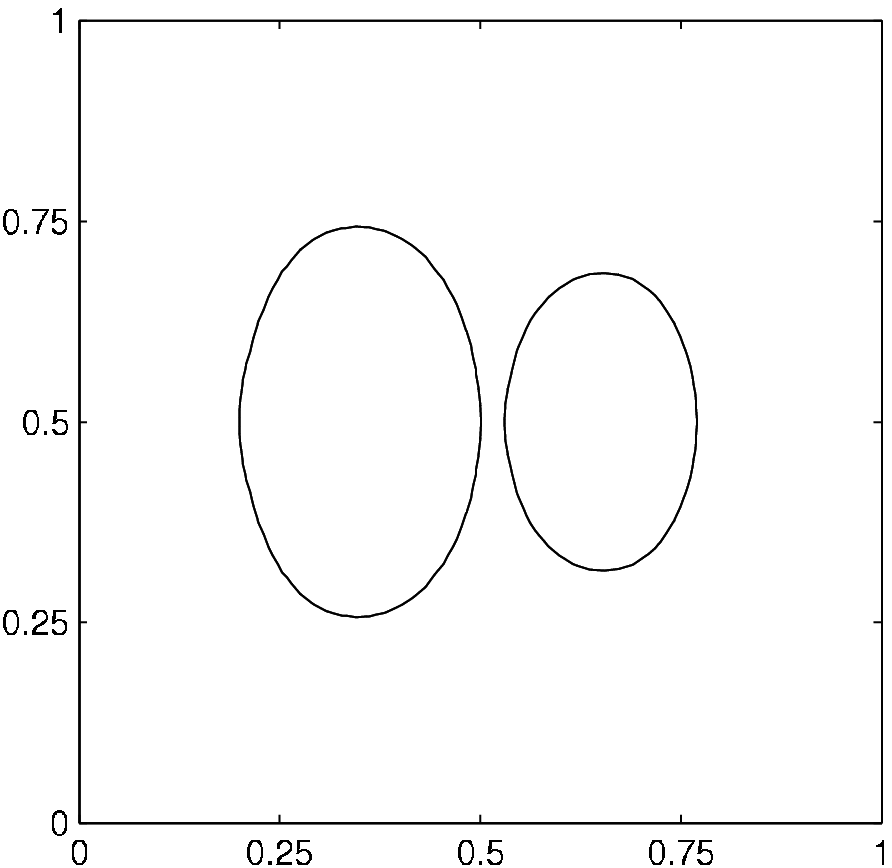} \\	
	\hspace{0.03cm}
	\includegraphics[width=3.9cm]{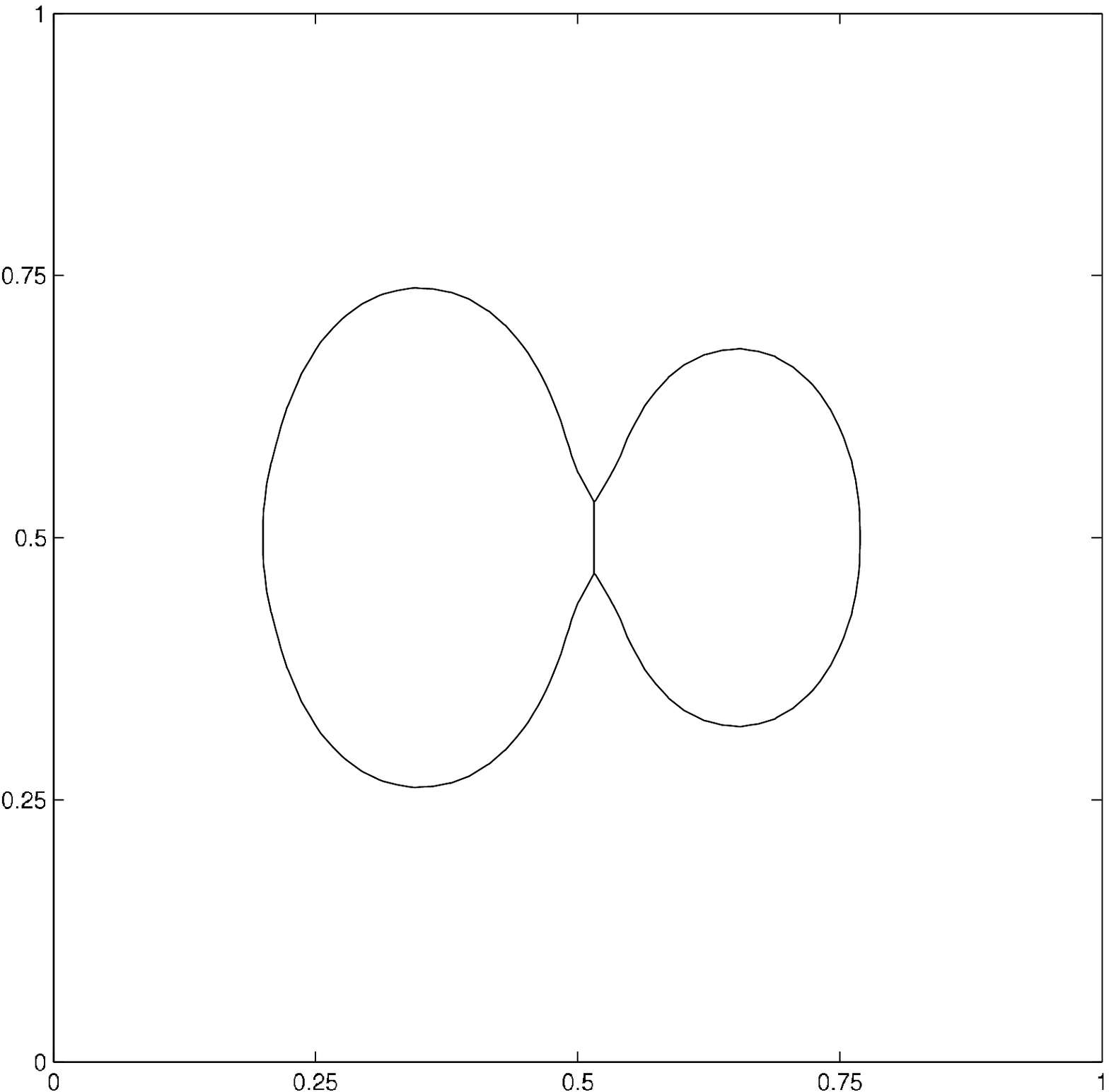} \hspace{1.02cm} 
	\includegraphics[width=4cm]{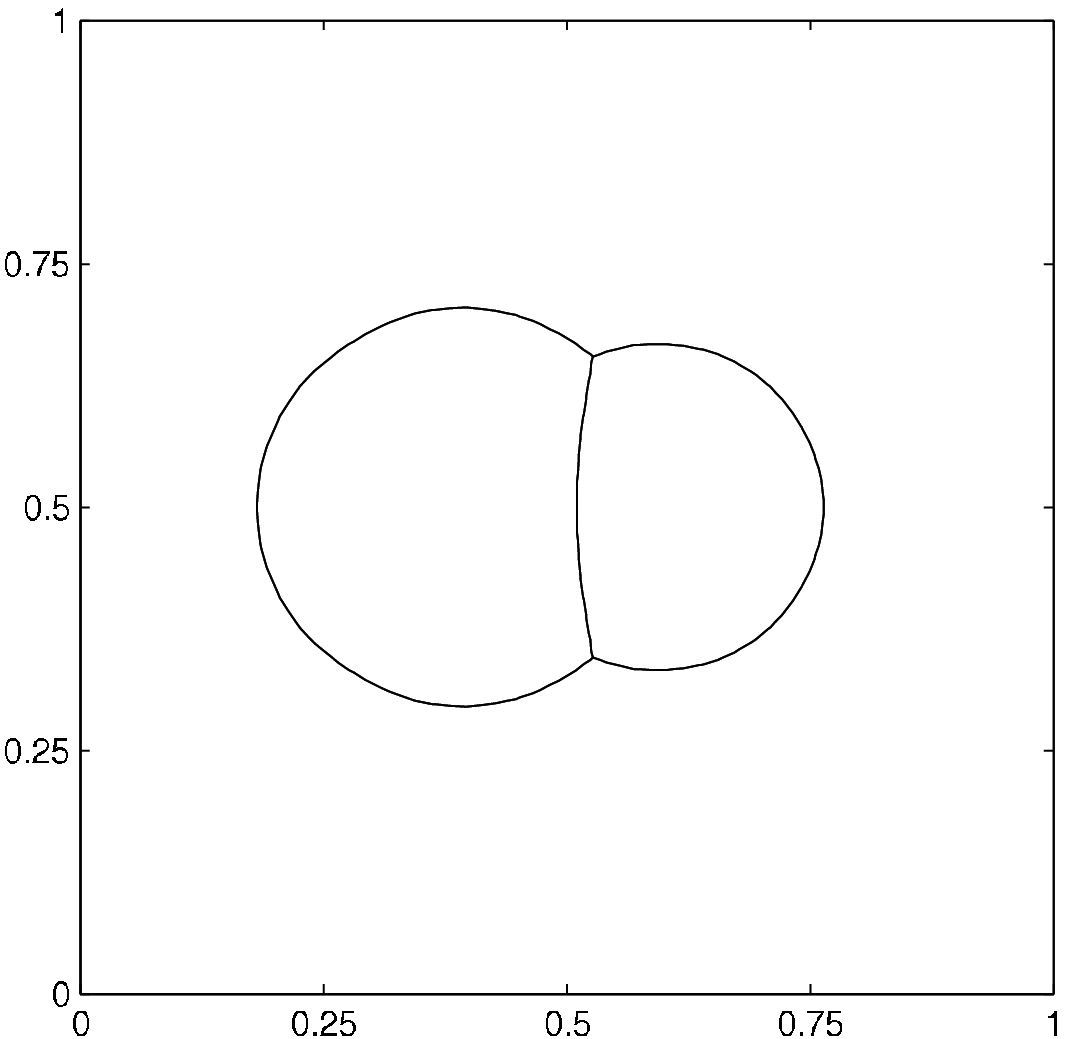}
	 	\caption{Attaching of two regions corresponding to different phases.}\label{fig:attach}
	 	\end{center}
\end{figure}

\subsection{An additional transport term}
Here we deal with a generalization of the constrained mean curvature flow to include a transport term.
For simplicity, we shall focus on the two-phase case which can be described by a scalar model.

Let us consider the partial differential equation (compare to (\ref{nonlinheat})),
\begin{align}
\label{eqtr}
	u_{t} &= \Delta u + \frac{f}{\sqrt{4 \pi t}} 
         +\lambda\mathcal{H}^{m-1}\lfloor_{\partial \{u> 1/2\}},
\end{align}
where $f$ is a specified function of space.
It can be shown in a fashion similar to \ref{appA} that the application of the BMO process to (\ref{eqtr}) leads to the motion of interfaces with normal velocity
\begin{align}
	v = - \kappa - f + \tilde{\lambda}, \notag
\end{align}
where $\kappa$ is the mean curvature and $\tilde{\lambda}$ is a function of time that guarantees the preservation of area.

In the numerical computations, one uses the method described in section \ref{numalg} minimizing the following penalized functional 
over $H^{1}(\Omega)$:
\begin{align}
	\mathcal{F}_{n}({w}) &= \int_{\Omega} \Big( \frac{ |{w}-{w}_{n-1}|^2}{2h}
		+\frac{| \nabla {w}|^2}{2}+\frac{fw}{\sqrt{4 \pi nh}} \Big) {dx} + \frac{1}{\epsilon} | A - A(w) |,\notag
\end{align}
where $A$ is the enclosed area of the bubble which should be preserved over time and $A(w)$ is the measure of the set $\{ {w} \ge 1/2\}$.

We apply the explained method to carry out a simple two-dimensional simulation
of gas bubbles rising from the bottom of a container filled with a viscous liquid (see also \cite{Ginder}).
In this case we set $f= \beta y$, where $y$ is the coordinate direction of gravity and $\beta$ is a constant expressing buoyancy.
We consider the case of a bubble having the shape of a partial ellipse and initially positioned at the bottom of the container.
Figure 11 shows the evolution at four distinct times for two different initial shapes.
The results were obtained using the parameters $\Delta t = 10^{-3}$, $K=20$, $\epsilon = 0.001$, and $\beta=20.5$. 

We compute under Neumann boundary conditions which means that the bubble
will always touch the bottom with right angle.
The motion is a result of the balance between the buoyant force pushing the bubble upwards and the surface tension force pressing the bubble towards the bottom.
For the bubble on the left, the adhesive force is prevalent so it attaches itself to the boundary and then comes to a rest in a stationary shape.
The bubble on the right has a shape for which the buoyant force wins and the bubble detaches itself from the boundary. 
After detaching, the bubble becomes circular and continues moving upward. 

\begin{center}
\begin{figure}[!ht]
\includegraphics*[
bb=81 227 529 564, 
trim=0 20 100 200, clip,
scale=0.85]{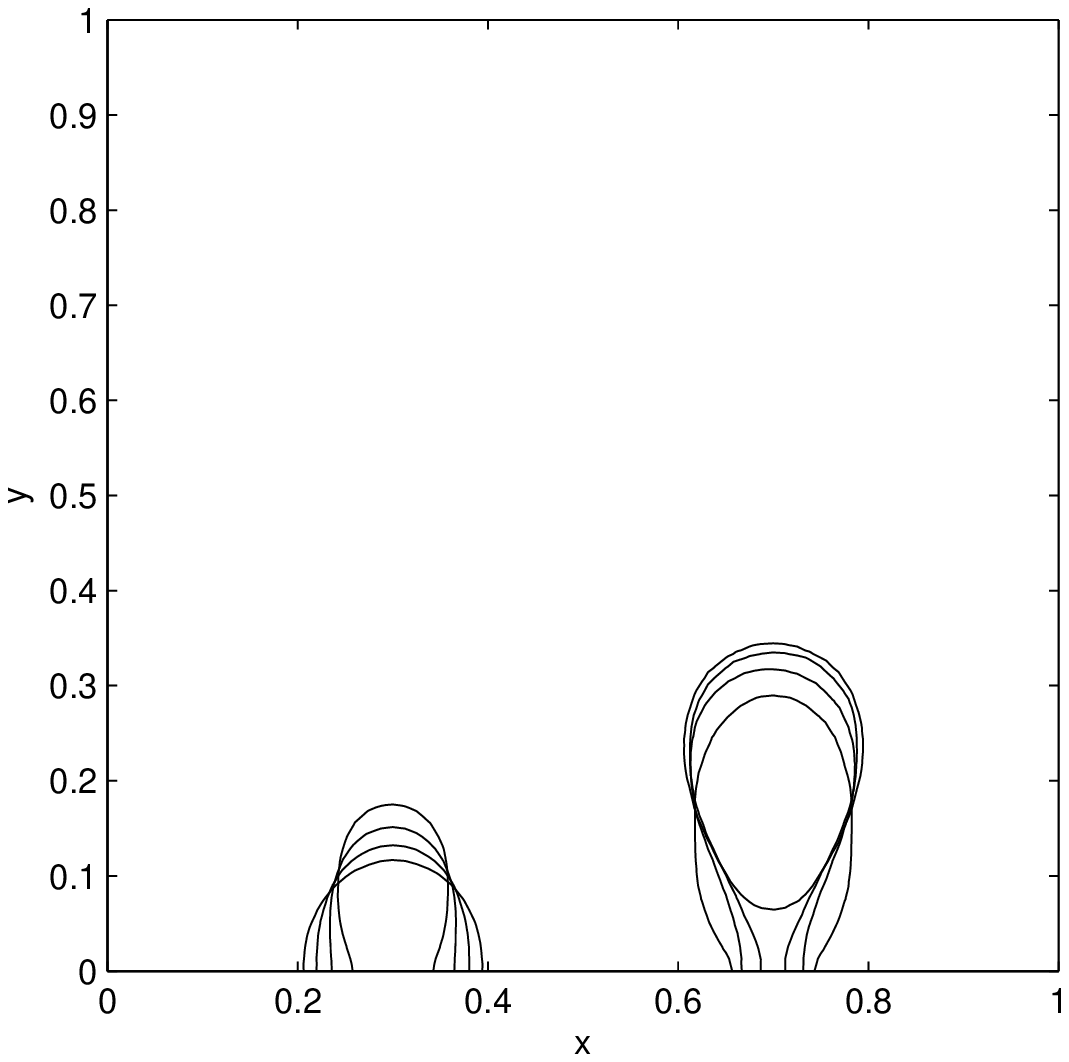}
   \caption{Volume-preserving mean curvature flow with buoyancy - comparison of evolution for different initial shapes.}
   \label{fig_nctam}
\end{figure}
\end{center}
\vspace{-10mm}

The possibility of a straightforward inclusion of a transport term into the numerical algorithm as explained above is extremely important in applications where two or more phases interact at the interface, such as through buoyancy in our simple example.
We believe that it is possible to proceed in this direction and design effective algorithms for various interaction problems.

\section{Conclusion}
We developed a method for approximating constrained multiphase cur\-va\-tu\-re-driven motions.
Our method is based on the BMO algorithm, which was reformulated in terms of a vector-valued heat equation.
We derived the nonlinear PDE which governs the corresponding constrained evolutions and used it to formally prove convergence in the multiphase setting.

The vector-valued BMO algorithm was implemented employing the discrete Morse flow and it was found that the variational nature of this approach allows one to consider constraints via additional penalty terms.
Using this idea we were able to realize multiphase area-preserving interfacial motions that are free of the defects of other methods.

By detecting the precise locations of the interfaces we are able to compute the area of each phase at a high precision and thus impose the area constraints. 
This geometric information is also retained after the thresholding step, which was found to alleviate the BMO's standard restrictions on the spatial and time-step resolutions, for both the constrained and non-constrained problems. 

In closing, we remark that it would be interesting to investigate our algorithm in relation to the recent threshold dynamics utilizing signed distances functions \cite{Esedoglu2}, and to consider its position in applications.

\appendix

\section{Construction of reference vectors}
\label{app_refvec}
For the sake of completeness, we include a method for constructing the vectors corresponding to the regular simplexes.

The computation of the reference vectors can be done by first considering the standard $(k-1)$-simplex in ${\mathbb R}^k$ and rotating its vertex vectors in a suitable way. Particularly, the standard $(k-1)$-simplex is given by
$$ S_{k-1} = \{ (x_1,x_2, \dots , x_k) \in {\mathbb R}^k; \;\; \sum_{i=1}^k x_i = 1, \; \text{and} \;  x_i \geq 0 \; \text{for} \; i=1, \dots , k \} . $$
Its vertices have the coordinates
$$ (1,0, \dots, 0), \; (0,1,0, \dots ,0), \dots , \; (0, \dots ,0,1) $$
and if we translate the simplex in such a way that its centroid lies in the origin, the vertices will be
\begin{eqnarray*}
\boldsymbol{p}^*_1&=& \tfrac{1}{k} (k-1,-1, \dots, -1), \\
\boldsymbol{p}^*_2 &=& \tfrac{1}{k}(-1,k-1,-1, \dots ,-1), \\
& \vdots & \\
\boldsymbol{p}^*_k &=& \tfrac{1}{k}(-1, \dots ,-1,k-1) .
\end{eqnarray*}
We fix an orthonormal basis for the $(k-1)$-dimensional hyperplane containing the simplex. A convenient way is to take the first translated vertex above as the first basis vector (after normalization) and construct the remaining vectors as follows:
\begin{eqnarray*}
\boldsymbol{q}_1^k &=& \tfrac{1}{\sqrt{(k-1)k}} (k-1,-1, -1,-1,\dots, -1) \\
\boldsymbol{q}_2^k &=& \tfrac{1}{\sqrt{(k-2)(k-1)}} (0,k-2,-1, -1,\dots, -1) \\
\boldsymbol{q}_3^k &=& \tfrac{1}{\sqrt{(k-3)(k-2)}} (0,0,k-3,-1, \dots, -1) \\
&& \dots \\
\boldsymbol{q}_{k-1}^k &=& \tfrac{1}{\sqrt{2}} (0, \dots, 0,1, -1) .
\end{eqnarray*}
Denoting $Q^k$ the matrix having $\boldsymbol{q}_1^k, \dots, \boldsymbol{q}_{k-1}^k$ as its rows, we obtain the reference vectors as the normalized projection of $\boldsymbol{p}^*_i,i=1, \dots,k$, into this orthonormal system, i.e, 
$$ \boldsymbol{p}_i^T = \frac{1}{| Q^k (\boldsymbol{p}^*_i)^T |} Q^k (\boldsymbol{p}^*_i)^T . $$

\begin{figure}[!ht]
\begin{center}
\includegraphics[height=4.1cm]{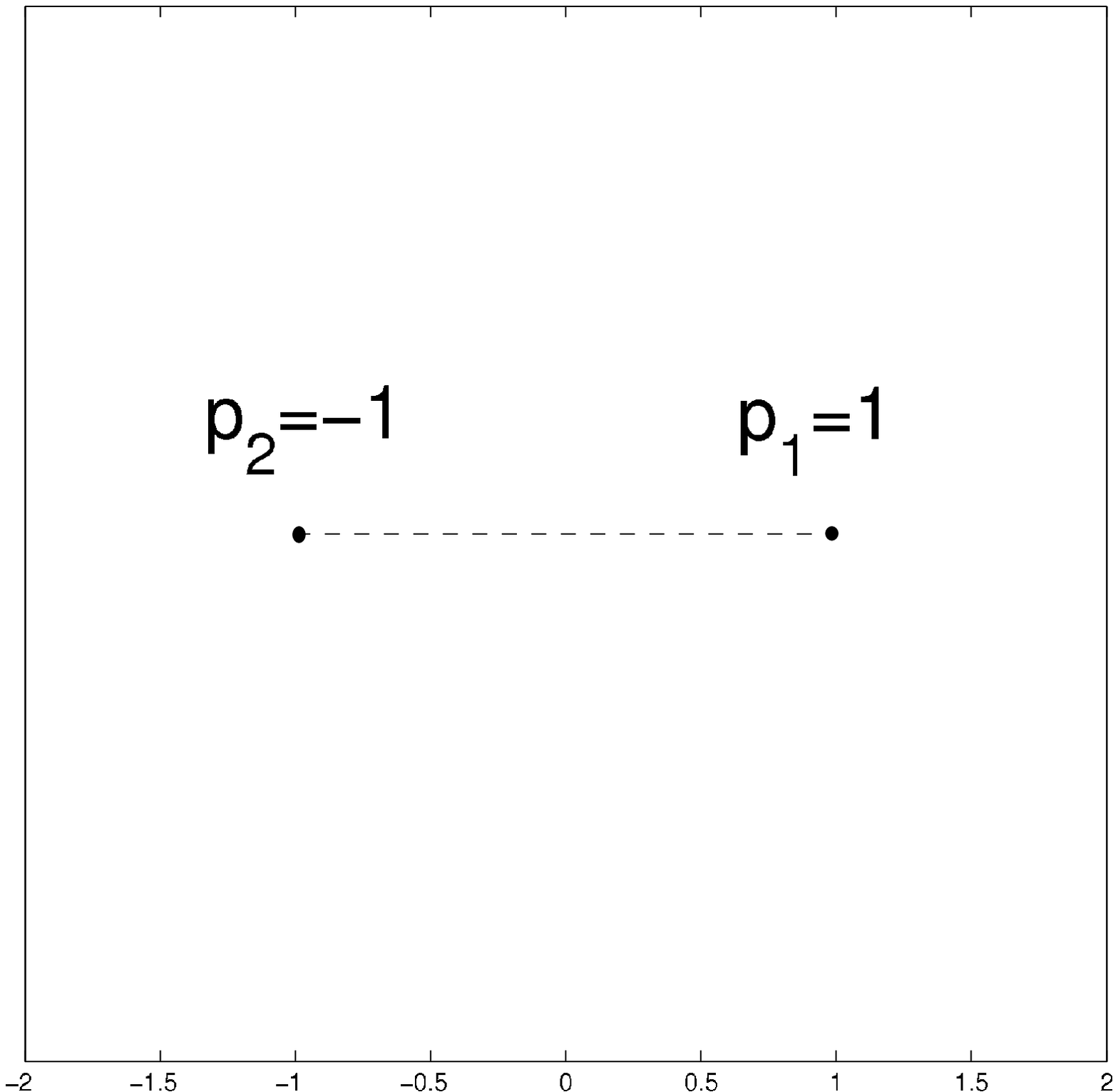} 
\hspace{0.2cm}
\includegraphics[height=4.1cm]{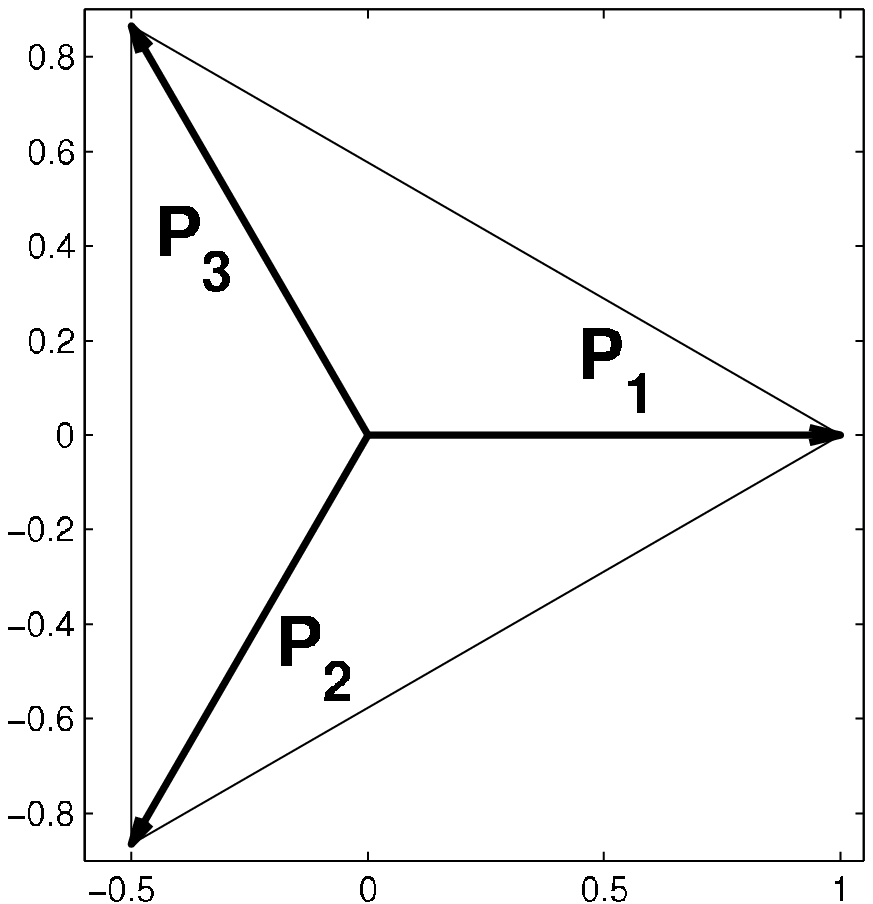}
\vspace{0.2cm}
\includegraphics[height=4.1cm]{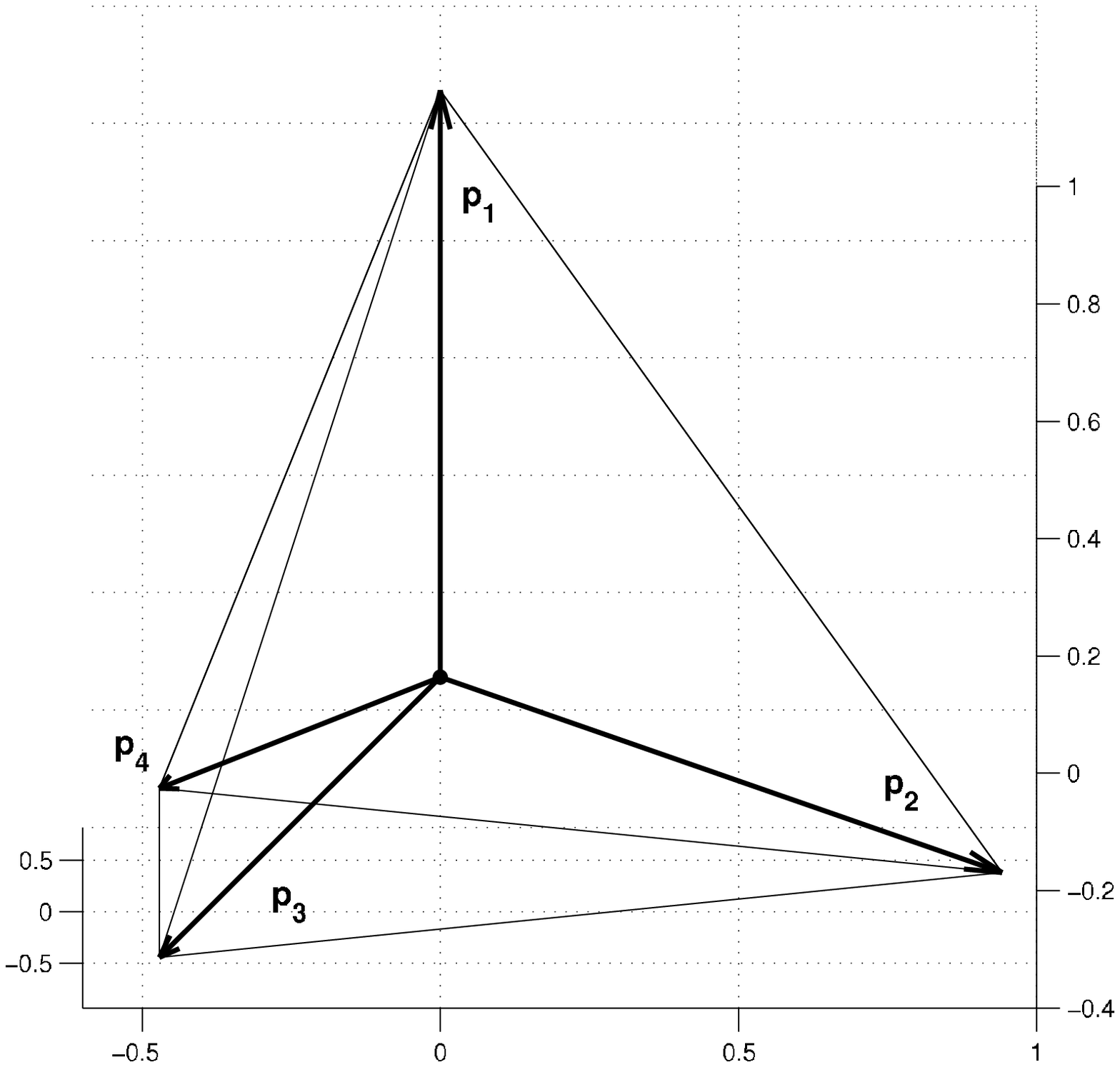}
\caption{(Left) The 2-phase regular simplex. (Center) The 3-phase regular simplex. (Right) The 4-phase regular simplex. }
\label{regsim}
\end{center}
\end{figure}

\section{Formal proof of convergence of the modified BMO algorithm}
\label{appA}

We show that the $\tfrac{1}{2}$-level set of the solution to the problem
\begin{eqnarray*}
u_t &=& \Delta u + \mu \qquad \; \text{in} \; (0,T) \times {\mathbb R}^2, \\
u(0,x) &=& \chi_{P(0)}(x) \qquad \text{in} \; {\mathbb R}^2,
\end{eqnarray*}
evolves according to its curvature $\kappa$, plus a constant factor $\lambda$, plus an error term which approaches $0$ as $T \to 0$.
Here $P(0)$ is a given initial region, $P(t)$ denotes $\{ x: \; u(x,t) > \tfrac{1}{2} \}$, and $\mu$ is a Radon measure given by
$$ \mu (t,x) = \lambda (t) {\mathcal H}^1\lfloor_{\partial P(t)} $$
for a suitable function $\lambda$. 
We assume that $\lambda(0)$ can be defined so that $\lambda$ is a smooth function in $[0,T]$.

\begin{figure}[!ht]
\label{fig_proof}
\begin{center}
\includegraphics[scale=0.7]{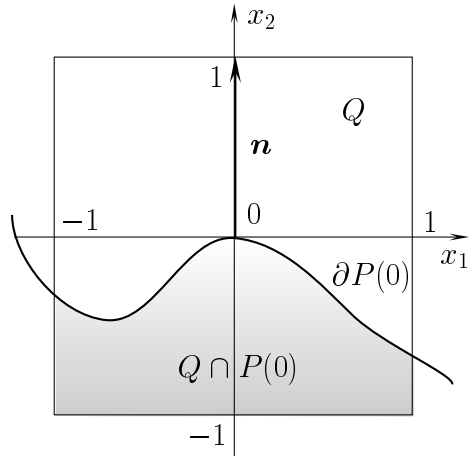}
\caption{Configuration of the interface in the proof of BMO convergence.}
\end{center}
\end{figure}

We consider the coordinate system as in figure B.13, where the point $(0,0)$ lies on $\partial P(0)$ and the outer normal $\boldsymbol{n}$ to $P(0)$ at this point is the vector $(0,1)$.
We assume that inside the cube $Q=\{ (x_1,x_2): \; |x_1| \leq 1, |x_2| \leq 1 \}$ the boundary of $P(0)$ is given by the graph of a function $\gamma: (-1,1) \to (-1,1)$, satisfying
$$ |\gamma'(x_1)| \leq 1, \; x_1 \in (-1,1), \qquad \gamma(0)=\gamma'(0)=0 . $$
Then $\gamma''(0)$ is equal to minus the curvature $\kappa$ of $\partial P(0)$ at the point $(0,0)$.
Moreover, we assume that for a sufficiently short period of time, $\partial P(t)$ can be written as the graph of a function $\gamma(t,x_1)$ in the same coordinate system.

Let $v$ be the normal velocity of $P(0)$, i.e., 
$$ vt  \boldsymbol{n} \in \partial P(t) \qquad \text{or} \qquad u(t,0,vt) = \tfrac{1}{2} .$$
Writing the explicit solution to the heat equation with a source term, we have that $u(t,0,vt)$ is equal to
$$ \frac{1}{4 \pi t} \int_{P(t)} e^{- \frac{x_1^2 + (x_2-vt)^2}{4t}} dx + \int_0^t \int_{{\mathbb R}^2} \frac{1}{4 \pi (t-s)} e^{- \frac{x_1^2 + (x_2-vt)^2}{4(t-s)}}\mu(s,x) \, dx \, ds . $$
By the results of \cite{Evans,Ishii}, the first integral on the right-hand side is equal to
$$ \frac{1}{2} + \frac{\sqrt{t}}{2 \sqrt{\pi}} (\gamma''(0)-v) + O(t^{3/2}) , \quad t \to 0+. $$
Let us denote the second integral on the right-hand side by $I$. Our goal is to show that
$$ I = \frac{\sqrt{t}}{\sqrt{\pi}} \lambda (0) + O(t^{3/2}) , \quad t \to 0+.$$

We split the integration domain into two parts as follows:
\begin{eqnarray*}
I &=& \int_0^t \int_{\partial P(s)} \frac{\lambda (s)}{4 \pi (t-s)} e^{- \frac{x_1^2 + (x_2-vt)^2}{4(t-s)}} dS(x) \, ds \\
&=& \int_0^t \int_{\partial P(s) \cap Q} \frac{\lambda (s)}{4 \pi (t-s)} e^{- \frac{x_1^2 + (x_2-vt)^2}{4(t-s)}} dS(x) \, ds \\
&& \qquad \qquad + \int_0^t \int_{\partial P(s) \setminus Q} \frac{\lambda (s)}{4 \pi (t-s)} e^{- \frac{x_1^2 + (x_2-vt)^2}{4(t-s)}} dS(x) \, ds  \\
&=& I_1 + I_2 .
\end{eqnarray*}

First we show that the second integral $I_2$ is exponentially small:
\begin{eqnarray*}
|I_2| &\leq& \int_0^t \int_{\partial P(s) \setminus Q} \frac{|\lambda (s)|}{4 \pi (t-s)} e^{\frac{-1}{4(t-s)}} dS(x) \, ds \\
&\leq& \max_{s \in [0,t]} | \partial P(s)| \int_0^t \frac{|\lambda (s)|}{4 \pi (t-s)} e^{\frac{-1}{4(t-s)}} ds \\
&=& C \int_{1/(4t)}^{\infty} | \lambda (t- \tfrac{1}{4 \sigma})| \frac{1}{\sigma} e^{- \sigma}  \, d \sigma \\
&\leq& 4Ct \max_{s \in [0,t]} |\lambda (s)| \int_{1/(4t)}^{\infty} e^{- \sigma} \, d \sigma \\
&=& Ct e^{-\frac{1}{4t}} \\
&=& O(e^{-\alpha/t}),
\end{eqnarray*}
for some $\alpha > 0$. 
In the above calculations we use the change of variable $\sigma = 1/4(t-s)$ and we assume that $\partial P(s)$ has finite length for $s \in [0,T]$ and that $\lambda (s)$ is bounded for $s \in [0,T]$.

In the case of $I_1$ we use the fact that the boundary of the domain is a graph and employ Green's theorem. 
We use the notation $\rho = \sqrt{t-s}$ in order to simplify the formulae.
Further, we perform a change of variables $ z_1 = x_1/ \rho$, $z_2 = (x_2- vt)/ \rho$ to obtain the value of $I_1$ as 
\begin{eqnarray*}
&& \int_0^t \int_{-1}^1 \int_{- \infty}^{\gamma (s,x_1)} \frac{\partial}{\partial x_2} \left[ \frac{\lambda (s)}{4 \pi (t-s)} e^{- \frac{x_1^2 + (x_2-vt)^2}{4(t-s)}} \right] \sqrt{1+ \gamma'(s,x_1)^2} \,  dx_2 \, dx_1 \, ds \\
&=& \int_0^t \int_{-1}^1 \int_{- \infty}^{\gamma (s,x_1)} \frac{\lambda (s)}{4 \pi \rho^2} e^{- \frac{x_1^2 + (x_2-vt)^2}{4(t-s)}} \frac{-x_2+vt}{2 \rho^2} \sqrt{1+ \gamma'(s,x_1)^2} \,  dx_2 \, dx_1 \, ds \\
&=& \int_0^t \int_{-1/\rho}^{1/\rho} \int_{- \infty}^{\frac{\gamma (s, \rho z_1)-vt}{\rho}}
\frac{\lambda (s)}{8 \pi \rho} e^{- \frac{z_1^4+z_2^2}{4}} (-z_2) 
\sqrt{1+ \gamma'(s,\rho z_1)^2} \,  dz_2 \, dz_1 \, ds \\
&=& \int_0^t \int_{-1/\rho}^{1/\rho} \int_{- \infty}^0 \frac{- \lambda (s) z_2}{8 \pi \rho} e^{- \frac{|z|^2}{4}} \sqrt{1+ \gamma'(s,\rho z_1)^2} \,  dz_2 \, dz_1 \, ds \\
&& \quad + \int_0^t \int_{-1/\rho}^{1/\rho} \int_0^{\frac{\gamma (s, \rho z_1)-vt}{\rho}} \frac{- \lambda (s) z_2}{8 \pi \rho} e^{- \frac{|z|^2}{4}} \sqrt{1+ \gamma'(s,\rho z_1)^2} \,  dz_2 \, dz_1 \, ds \\
&=& I_{11} + I_{12} .
\end{eqnarray*}

In the sequel, we shall need the asymptotic properties of the arc-length term $\sqrt{1+\gamma'(s,\rho z_1)^2}$. Taylor's expansion at $s=0, z=0$ gives
\begin{eqnarray*}
&& \sqrt{1+\gamma'(s,z)^2} = \sqrt{1+\gamma'(0,0)^2} + \frac{\gamma'(0,0)\gamma'_t(0,0)}{\sqrt{1+\gamma'(0,0)^2}} s + \frac{\gamma'(0,0)\gamma''(0,0)}{\sqrt{1+\gamma'(0,0)^2}} z  \\
&& \qquad \qquad \qquad \qquad \qquad \qquad + F_1(\tau,\xi)s^2 + F_2(\tau,\xi)sz + F_3(\tau,\xi)z^2, 
\end{eqnarray*}
where $F_1,F_2,F_3$ are functions that are assumed to be smooth and bounded on the segment $(\tau,\xi)$ connecting the points $(0,0)$ and $(s,z)$.
Hence we can write
$$ \sqrt{1+\gamma'(s,z_1 \sqrt{t-s})^2} = 1 + O(s^2+(t-s)z_1^2) . $$
In the same way,
$$ \gamma(s,z) = \gamma(0,0) + \gamma_t(0,0)s + \gamma'(0,0)z+ O(s^2 + z^2) = vs + O(s^2+z^2) , $$
i.e.,
$$ \gamma(s,z_1 \sqrt{t-s}) = vs + O(s^2+(t-s)z_1^2). $$

Let us denote for simplicity
$$ \eta := \frac{\gamma(s, \rho z_1) - vt}{\rho} =  - v \sqrt{t-s} + O( \frac{s^2}{\sqrt{t-s}} + z_1^2 \sqrt{t-s}) .$$
Again, the second integral $I_{12}$ is small. Since $exp(-z_2^2) \leq 1$, we have
\begin{eqnarray*}
|I_{12}| &\leq& C \int_0^t \int_{-1/\rho}^{1/\rho} \int_0^{\eta} \frac{| \lambda (s)|}{\sqrt{t-s}} |z_2| e^{- \frac{z_1^2}{4}} \big(1+ O(s^2+(t-s)z_1^2) \big) \,  dz_2 \, dz_1 \, ds \\
&\leq& C \int_0^t \int_{-1/\rho}^{1/\rho} \frac{| \lambda (s)|}{\sqrt{t-s}} \eta^2 e^{- \frac{z_1^2}{4}} \big(1+ s^2+(t-s)z_1^2 \big) \, dz_1 \, ds \\
&\leq& C \int_0^t \int_{-\infty}^{\infty} \frac{| \lambda (s)|}{\sqrt{t-s}} \big( (t-s)(v^2+z_1^4) + \frac{s^4}{t-s}  \big) \\
&& \qquad \qquad \qquad \qquad \qquad \cdot e^{- \frac{z_1^2}{4}} \big(1+ s^2+(t-s)z_1^2 \big) \, dz_1 \, ds \\
&\leq& C \int_0^t \frac{| \lambda (s)|}{\sqrt{t-s}} \big( (t-s) + \frac{s^4}{t-s} \big) \big(1+ s^2+(t-s) \big) \, ds \\
&=& O(t^{3/2}) .
\end{eqnarray*}

It remains to compute $I_{11}$. Using the Taylor expansion we obtain
\begin{eqnarray*}
I_{11} &=& \int_0^t \int_{-1/\rho}^{1/\rho}  \frac{ \lambda (s)}{8 \pi \rho} e^{- \frac{z_1^2}{4}}\sqrt{1+ \gamma'(\rho z_1)^2} \Big( \int_{- \infty}^0 (-z_2) e^{- \frac{z_2^2}{4}} \,  dz_2 \Big) \, dz_1 \, ds \\
&=& \int_0^t \int_{-1/\rho}^{1/\rho}  \frac{ \lambda (s)}{4 \pi \rho} e^{- \frac{z_1^2}{4}} (1+ O(s^2+(t-s)z_1^2) ) \, dz_1 \, ds \\
&=& \int_0^t \frac{ \lambda (s)}{4 \pi \sqrt{t-s}} \int_{-\infty}^{\infty} e^{- \frac{z_1^2}{4}} \, dz_1 \, ds  - \int_0^t \frac{ \lambda (s)}{4 \pi \sqrt{t-s}} \int_{-\infty}^{-1/\rho} e^{- \frac{z_1^2}{4}} \, dz_1 \, ds \\
&& \qquad - \int_0^t \frac{ \lambda (s)}{4 \pi \sqrt{t-s}} \int_{1/\rho}^{\infty} e^{- \frac{z_1^2}{4}} \, dz_1 \, ds \\
&& \qquad + \int_0^t \frac{ \lambda (s)}{4 \pi \sqrt{t-s}} \int_{-1/\rho}^{1/\rho}  O(s^2+(t-s)z_1^2) e^{- \frac{z_1^2}{4}} \, dz_1 \, ds \\
&=& I_{111} + I_{112} + I_{113} + I_{114}
\end{eqnarray*}

It turns out that the first integral $I_{111}$ gives the desired value, while the remaining integrals are small. To see this, we write
\begin{eqnarray*}
I_{111} &=& \int_0^t \frac{ \lambda (s)}{\sqrt{4 \pi} \sqrt{t-s}}  \, ds 
= \int_0^t \frac{\lambda(0)+ O(s)}{\sqrt{4 \pi} \sqrt{t-s}} \, ds 
= \frac{\sqrt{t}}{\sqrt{\pi}} \lambda (0) + O(t^{3/2}), \\
|I_{112}| &\leq& \int_0^t \frac{| \lambda (s)|}{4 \pi} \int_{-\infty}^{-1/\rho} |z_1| e^{- \frac{z_1^2}{4}} \, dz_1 \, ds 
= \int_0^t \frac{| \lambda (s)|}{2 \pi} e^{- \frac{1}{4(t-s)}} \, ds 
= O(e^{-\alpha/t}).
\end{eqnarray*}
The estimate for $I_{113}$ is similar and
\begin{eqnarray*}
|I_{114}| &\leq& C \int_0^t \frac{ |\lambda (s)|}{ \sqrt{t-s}} \int_{-\infty}^{\infty} (s^2+(t-s)z_1^2) e^{- \frac{z_1^2}{4}} \, dz_1 \, ds \\
&\leq& C \int_0^t \frac{|\lambda (s)| }{\sqrt{t-s}}(s^2+t-s) \, ds \\
&=& O(t^{3/2}).
\end{eqnarray*}

Gathering the results, from the relation $u(t,0,vt)=1/2$ we have arrived at the identity
$$ \frac{1}{2} + \frac{\sqrt{t}}{2 \sqrt{\pi}} (- \kappa - v) + \frac{\sqrt{t}}{\sqrt{\pi}} \lambda (0) + O(t^{3/2}) = \frac{1}{2}, \qquad t \to 0+. $$
Therefore,
$$ v = -\kappa + 2 \lambda (0) + O(t), \qquad t \to 0+.$$

\section{Derivation of conditions in the multiphase setting}
\label{app_multi}

We derive the equation and free boundary conditions for the case of multiple phases.
Since the technical aspect is similar to that of the two-phase case, we will give only the main ideas.

Let us denote the phase regions by $P_i$, $i=1, \dots , k$, and the interface between regions $P_i$ and $P_j$ by $\gamma_{ij}$.
The symbol $u_i$ will mean the value of the inner product $\boldsymbol{u} \cdot \boldsymbol{p}_i$, where $\boldsymbol{p}_i$ is one of the reference vectors constructed in Section \ref{multinoBMO}.  
The simplest situation for three phases is depicted in figure C.14.
Here the whole sets $\{u_i = u_j \}$, $i,j=1,2,3$ are shown, and the parts of these sets corresponding to the interfaces are drawn with thicker lines.

\begin{figure}[!ht]
\label{fig_fb_multi}
\begin{center}
\includegraphics*[
bb=71 525 697 721, 
trim=100 0 205 0, clip,
scale=0.72]{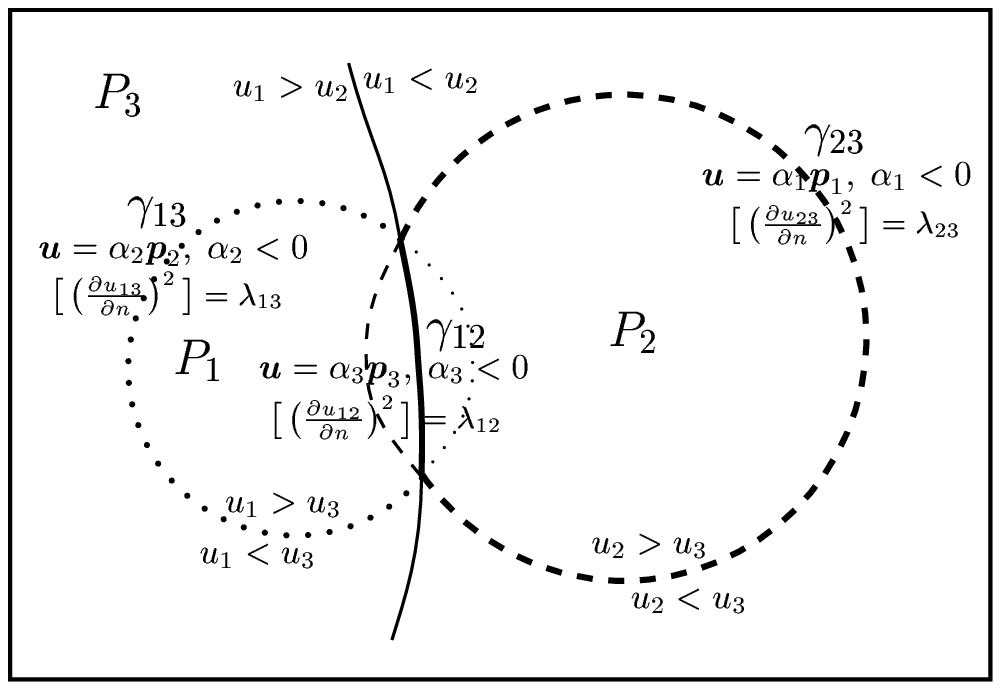}
\includegraphics*[scale=0.5]{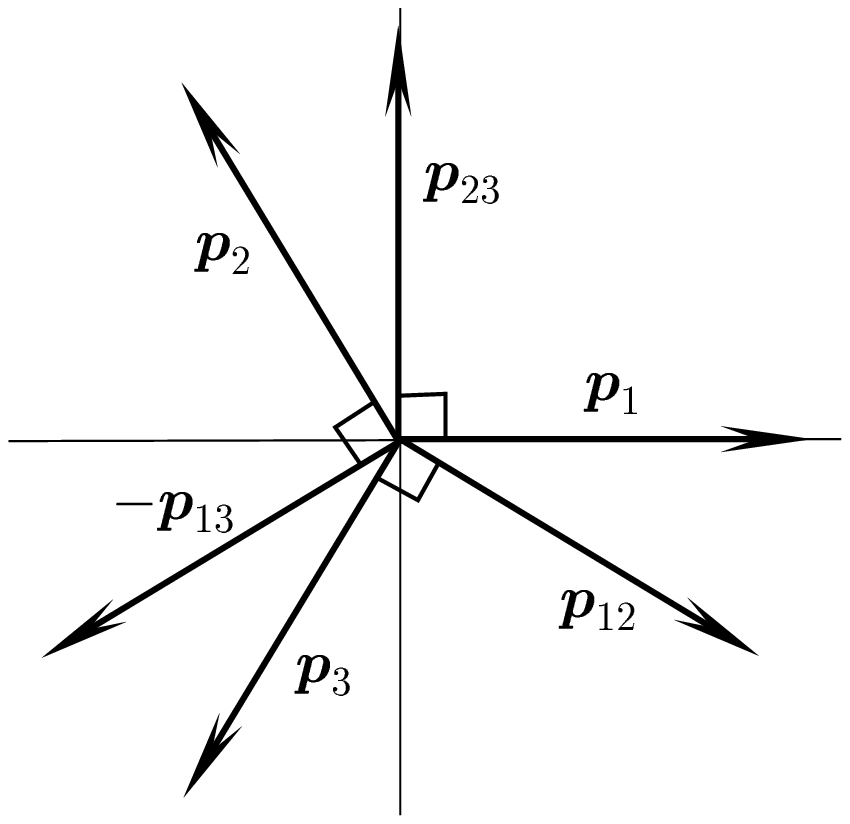}
\caption{Conditions holding on the interfaces in the case of 3 phases (left) and configuration of reference vectors (right).}
\end{center}
\end{figure}

We consider an arbitrary interface $\gamma_{ij}$.
A function $\boldsymbol{u}: \mathbb{R}^m \to \mathbb{R}^{k-1}$ ($m=2$ in the above figure) can be expressed as a linear combination of reference vectors in the following way:
\begin{equation}
\label{udec}
\boldsymbol{u}(x) = \sum_{l \neq i,j} \alpha_l(x) \boldsymbol{p}_l + \beta(x) \boldsymbol{p}_{ij} . 
\end{equation}
Here $\boldsymbol{p}_{ij}=(\boldsymbol{p}_i-\boldsymbol{p}_j)/|\boldsymbol{p}_i-\boldsymbol{p}_j|$ is a unit vector orthogonal to $\mathcal{L}_{ij} = span \{ \boldsymbol{p}_l , l \neq i,j \}$, and $\beta = u_{ij}$ where $ u_{ij} := \boldsymbol{u} \cdot \boldsymbol{p}_{ij}$.
Since $\boldsymbol{p}_{ij}=-\boldsymbol{p}_{ji}$ and $u_{ij}=-u_{ji}$, we will use the symbols $u_{ij}$ and $\boldsymbol{p}_{ij}$ only for indices $i,j$ satisfying $i<j$ to avoid confusion.

From the condition $u_i=u_j$ on $\gamma_{ij}$ and since $u_i>u_l$ inside $P_i$ whenever $i \neq l$, we obtain 
\begin{equation}
\label{conmul}
\beta = 0, \qquad \alpha_l <0 \;\;\; \forall l \neq i,j \qquad \text{on} \;\; \gamma_{ij}.
\end{equation}
On the other hand, using the relation $\boldsymbol{p}_i \cdot \boldsymbol{p}_j = -1/(k-1)$, we compute the Dirichlet functional as
\begin{eqnarray}
\label{dirmod}
J(\boldsymbol{u}) &=& \int_{\Omega} | \nabla \boldsymbol{u} |^2 \, dx \\
&=&  \int_{\Omega} \Big\{ \sum_{l \neq i,j} | \nabla \alpha_l |^2 - \frac{2}{k-1} \sum_{l,m \neq i,j} \nabla \alpha_l \cdot \nabla \alpha_m + | \nabla \beta |^2 \Big\} \, dx. \nonumber
\end{eqnarray}
In view of the separateness of conditions on $\alpha_l, \beta$ in (\ref{conmul}), we can proceed similarly as in the two-phase case to arrive at a free boundary problem corresponding to the steepest descent of $J$ under area constraints $|P_l| = A_l$, $l=1, \dots ,k$:
\begin{eqnarray}
\label{fbp1}
\boldsymbol{u}_t &=& \Delta \boldsymbol{u} \qquad \text{in} \;\; P_l, \; l=1, \dots, k, \\
\label{fbp2}
\frac{\partial \boldsymbol{u}}{\partial n} &=& 0 \qquad \quad \text{on} \;\; \partial \Omega , \\
\label{fbp3}
\left[ \Big( \frac{\partial u_{ij}}{\partial n} \Big)^2 \right]_{\gamma_{ij}} &=& \lambda_{ij} \qquad  \; i,j = 1, \dots, k, \;\; i<j, \\
\label{fbp4}
\left[ \frac{\partial u_l}{\partial n} \right]_{\gamma_{ij}} &=& 0 \qquad \;\;\;\; l \neq i,j, \;\; i,j =1, \dots , k.
\end{eqnarray}
Here $[ w ]_{\gamma_{ij}}$ denotes the jump of $w$ across $\gamma_{ij}$ in the normal direction and $\lambda_{ij}$'s are suitable functions of time.

We assert that the above derived free boundary problem arises from the unconstrained gradient flow of the functional
$$ J^{\tilde{\lambda}} ( \boldsymbol{u}) = \int_{\Omega} \Big\{ | \nabla \boldsymbol{u} |^2 + \sum_{l=1}^{k-1} \tilde{\lambda}_l \prod_{m \neq l}\chi_{\boldsymbol{u} \cdot \boldsymbol{p}_l > \boldsymbol{u} \cdot \boldsymbol{p}_m} \Big\} \, dx  .$$
This assertion is natural from a formal standpoint of the theory of Lagrange multipliers, since only the phase areas, multiplied by Lagrange multipliers $\tilde{\lambda}_l$, are added to the Dirichlet integral.
However, this can be also proved by calculations of the first and inner variations of the functional.

We present just an outline of the calculation of (\ref{fbp3}) for a selected interface $\gamma_{ij}$, $i<j<k$. 
Equations (\ref{fbp1}) and (\ref{fbp2}) are immediate and the identity (\ref{fbp4}) is derived in an analogous way as (\ref{fbp3}).
The idea is to decompose $\boldsymbol{u}$ again as in (\ref{udec}) and consider the following perturbation:
$$ \boldsymbol{u}^{\varepsilon} (x) = \sum_{l \neq i,j} \alpha_l(x) \boldsymbol{p}_l + u_{ij}(\eta^{-1}_{\varepsilon}(x)) \boldsymbol{p}_{ij}, $$
where $\eta_{\varepsilon}: {\bf R}^m \to {\bf R}^m$ is a function defined by
$$ \eta_{\varepsilon} (x) = x + \varepsilon \zeta (x) , $$
with $\zeta : {\bf R}^m \to {\bf R}^m$ a smooth function supported in a neighborhood of the interface $\gamma_{ij}$ and a positive distance away from every other set $\{ u_l = u_m \}$.
Because of the location of the support of $\zeta$, all characteristic functions $\chi_{\boldsymbol{u} \cdot \boldsymbol{p}_l > \boldsymbol{u} \cdot \boldsymbol{p}_m}$ in the expression for $J^{\tilde{\lambda}}$ remain unaffected by the introduced perturbation except for the terms $\chi_{\boldsymbol{u} \cdot \boldsymbol{p}_i > \boldsymbol{u} \cdot \boldsymbol{p}_j}$ and $\chi_{\boldsymbol{u} \cdot \boldsymbol{p}_j > \boldsymbol{u} \cdot \boldsymbol{p}_i}$.
Therefore, using the reformulation of the Dirichlet integral (\ref{dirmod}), we have
\begin{eqnarray*}
J^{\tilde{\lambda}}( \boldsymbol{u}) - J^{\tilde{\lambda}}(\boldsymbol{u}^{\varepsilon}) &=&
\int_{\Omega} \Big\{ | \nabla u_{ij}|^2 - | \nabla u^{\varepsilon}_{ij} |^2 + \\
&& \quad 
\sum_{l=i,j} \tilde{\lambda}_l \prod_{m \neq l} \chi_{\boldsymbol{u} \cdot \boldsymbol{p}_l > \boldsymbol{u} \cdot \boldsymbol{p}_m} -
\sum_{l=i,j} \tilde{\lambda}_l \prod_{m \neq l} \chi_{\boldsymbol{u}^{\varepsilon} \cdot \boldsymbol{p}_l > \boldsymbol{u}^{\varepsilon} \cdot \boldsymbol{p}_m} \Big\} \, dx . 
\end{eqnarray*}
After a long technical computation, which we omit, we obtain a condition holding on the interface
\begin{eqnarray*}
\int_{\gamma_{ij}} \Big\{ \Big[ | \nabla u_{ij}|^2 \zeta - 2(\nabla u_{ij} \cdot \zeta) \nabla u_{ij} \Big]_{\gamma_{ij}} \cdot \nu + (\tilde{\lambda}_i - \tilde{\lambda}_j) (\zeta \cdot \nu) \Big\} \, dl =0, 
\end{eqnarray*}
where $\nu = - (\nabla u_{ij})_{P_i} / | (\nabla u_{ij})_{P_i}|$ is the unit outer normal to phase region $P_i$ on $\gamma_{ij}$.
Since $\zeta$ was arbitrary, this yields the identity (\ref{fbp3}) with $\lambda_{ij} = \tilde{\lambda}_i - \tilde{\lambda}_j$.

From the above results one can deduce that the nonlinear PDE corresponding to (\ref{nonlinheat}) in the multiphase setting will be
\begin{equation}
\label{nonlpdev}
\boldsymbol{u}_t = \Delta \boldsymbol{u} +  \sum_{i=1}^{k-1} \lambda_i \boldsymbol{p}_i {\mathcal H}^{m-1} \lfloor_{ \partial P_i}, 
\end{equation}
where the $\lambda_i$'s are piecewise constant on the interfaces:
$$ \lambda_i = \sqrt{\tfrac{k-1}{2k}} \sum_{j \neq i} (\tilde{\lambda}_i - \tilde{\lambda}_j) \prod_{l \neq i,j} \chi_{u_l <u_i} \chi_{u_l < u_j}. $$
Indeed, taking the inner product of equation (\ref{nonlpdev}) with the vector $\boldsymbol{p}_{lm}$, we use the orthogonality properties and find that
$$ (u_{lm})_t = \Delta u_{lm} + (\tilde{\lambda}_l - \tilde{\lambda}_m) {\mathcal H}^{m-1} \lfloor_{\gamma_{lm}} $$
in a neighborhood of $\gamma_{lm}$.
This means that the function $u_{lm}$ satisfies a scalar equation of the type (\ref{nonlinheat}) and thus the condition (\ref{fbp3}) holds for all admissible $i,j$.
In the same vein, we can use the calculation from \ref{appA} to show that each interface moves with normal velocity equal to minus mean curvature plus a space-independent term, resulting in area preservation. 

\bibliography{GOS_references}
\bibliographystyle{model1a-num-names}

\end{document}